\documentclass[11pt,technote,onecolumn]{IEEEtran}

\usepackage{amssymb}
\usepackage{hyperref}
\usepackage{amssymb}
\usepackage{bm}
\usepackage{subfigure}

\usepackage{multirow}
\usepackage{amsmath}
\usepackage[labelfont=bf]{caption}
\usepackage{graphicx}
\usepackage{array}

\usepackage{cite,amsmath,graphicx}
\usepackage{bm}
\usepackage{amssymb}
\usepackage{booktabs}
\usepackage{multirow}
\usepackage{epstopdf}
\usepackage{amsthm}
\usepackage{lineno}
\usepackage{mathrsfs}
\usepackage{amsmath}
\usepackage{txfonts}
\usepackage{caption}
\usepackage{url}

\usepackage{amsfonts}

\newtheorem{definition}{Definition}
\newtheorem{lemma}{Lemma}

\DeclareMathOperator*{\arcos}{arcos}

\ifCLASSINFOpdf
 
\else

\fi

\usepackage{algorithm,algorithmic}

\begin{document}
%
\title{Hierarchical Tensor Ring Completion\\}
%
%
%

\author{Abdul Ahad,
        Zhen~Long,
        Ce~Zhu,~\IEEEmembership{Fellow,~IEEE}
        and~Yipeng~Liu,~\IEEEmembership{Senior Member,~IEEE}
\thanks{This research is supported by National Natural Science Foundation of China (NSFC, No. 61602091, No. 61571102) and Sichuan Science and Technology Program (No. 2019YFH0008, No. 2018JY0035). The corresponding author is Yipeng Liu.}
\thanks{All the authors are with School of Information and Communication Engineering, University of Electronic Science and Technology of China (UESTC), Chengdu, 611731, China. email: yipengliu@uestc.edu.cn.}
}

\markboth{Journal Name,~Vol.~, No.~, Month~Year}%
{Shell \MakeLowercase{\textit{et al.}}: Bare Demo of IEEEtran.cls for IEEE Journals}

\maketitle

\begin{abstract}

Tensor completion can estimate missing values of a high-order data from its partially observed entries.  Recent works show that low rank tensor ring approximation is one of the most powerful tools to solve tensor completion problem. However, existing algorithms need predefined tensor ring rank which may be hard to determine in practice. To address the issue, we propose a hierarchical tensor ring decomposition for more compact representation. We use the standard tensor ring to decompose a tensor into several 3-order sub-tensors in the first layer, and  each sub-tensor is further factorized by tensor singular value decomposition (t-SVD) in the second layer. In the low rank tensor completion based on the proposed decomposition, the zero elements in the 3-order core tensor are pruned in the second layer, which helps to automatically determinate the tensor ring rank. To further enhance the recovery performance, we use total variation to exploit the locally piece-wise smoothness data structure. The alternating direction method of multiplier can divide the optimization model into several subproblems, and each one can be solved efficiently. Numerical experiments on color images and hyperspectral images demonstrate that the proposed algorithm outperforms state-of-the-arts ones in terms of recovery accuracy.

\end{abstract}


%
\IEEEpeerreviewmaketitle

\section{Introduction}
Tensor, which is the higher-order generalization of vector and matrix, provides a natural form to represent higher-order data.  For example, a color image has three indices that can be represented by a 3-order tensor.  The tensor for  multidimensional data processing has attracted much attention in different fields such as signal and image processing~\cite{momeni2019generalized, madathil2019tensor}, computer vision~\cite{lu2019low}, quantum chemistry~\cite{mutlu2019toward,li2019stress} and data mining~\cite{papalexakis2016tensors,kong2019federated}. However, in some applications, a part of entries of multidimensional data are missing during data acquisition or transmission, which has a significant influence on subsequent processing. Tensor completion can recover missing entries from its observed entries by exploit the coherence in multi-linear space. The  low-rank method is one of the most powerful ones for tensor completion problems~\cite{liu2019image,leurgans1993decomposition,acar2011scalable,andersson1998improving,liu2013tensor,gandy2011tensor,liu2019low,lu2019low,long2019low,huang2020provable}.

The low-rank tensor completion methods are mainly divided into two categories according to different tensor rank formats.  The first one is based on tensor factorization. In this group, the tensor rank is predefined and the goal is to optimize the factors of tensor decomposition. For example, in~\cite{leurgans1993decomposition} and ~\cite{acar2011scalable}, the missing values of data are recovered with given CP rank in advance and each factor is updated by alternating least square (ALS) and gradient methods. Following it, Tucker-ALS~\cite{andersson1998improving}  recovers the missing entries with known Tucker rank.  The second one is to directly minimize the tensor rank. However, the rank is a non-convex function and solving the rank problem is NP-hard. Most existing methods use different nuclear norms as convex surrogates to solve this problem. For instance, HaLRTC~\cite{liu2013tensor} has been proposed to recover missing elements, which directly minimizes Tucker rank, and uses the nuclear norm of the unfolding matrix to solve the Tucker rank problem.  Besides, a series of works~\cite{gandy2011tensor, tan2014tensor} are developed to achieve a better recovery performance with Tucker decomposition. Apart from CP rank and Tucker rank in the two groups \cite{gandy2011tensor, romera2013new, liu2013tensor, kressner2014low,  xu2015parallel, yang2016iterative, yang2016rank, signoretto2011tensor,mu2014square, zhao2015bayesian,yang2015rank}, some other tensor ranks are used too, such as tensor train rank~\cite{oseledets2011tensor, 7859390}, tensor tree rank~\cite{hackbusch2009new, ballani2013black, da2015optimization, rauhut2015tensor}, tensor ring rank~\cite{zhao2016tensor, huang2020provable, wang2017efficient} and tubal rank~\cite{kilmer2013third, zhang2014novel,  zhang2016exact}.

Among all the tensor decompositions, tensor networks can capture more correlations than the rest ones. The corresponding tensor network based completion methods show superior performance, such as STTC~\cite{liu2019image}, TMac-TT~\cite{bengua2017efficient} and TR-ALS~\cite{wang2017efficient}. However, tensor train rank has its entries large for the middle factors and small for border factors, which leads to an unbalanced decomposition. To alleviate the drawbacks of tensor-train, the generalization of tensor-train decomposition named tensor-ring decomposition has been proposed in \cite{zhao2016tensor}. Low-rank tensor-ring completion is a powerful tool to recover missing data.  In~\cite{wang2017efficient}, the authors directly optimize  tensor ring factors with predefined tensor ring rank. Following it, 
in~\cite{he2019remote}, the authors add total variation (TV) regularization to model the local structure of image for low rank tensor ring completion on remote sensing image reconstruction. However, these two algorithms need predefined tensor ring rank which is hard to determine in practice.  In~\cite{yuan2019tensor}, the authors propose another low rank tensor ring completion method which applies rank minimization regularization on 3-order factors. It results in high computational cost when data are in large-scale. 

In this paper, to improve the tensor-ring based works, we propose a hierarchical low-rank tensor ring decomposition. For the first layer, we use the traditional tensor-ring decomposition model to factorize a tensor into many 3-order sub-tensors. For the second layer, each 3-order tensor is further decomposed by the tensor singular value decomposition (t-SVD)~\cite{kilmer2013third}.  For each step, we prune the zero elements in the 3-order core tensor, which achieves automatically rank determination.
We use this advanced tensor network for exploiting the global multi-linear data structure in low rank tensor completion. In the proposed optimization model for tensor completion, we additional employ total variation to exploit the local similarity of data in the form of piecewise smoothness. The alternating direction method of multipliers (ADMM) is used to solve the optimization problem. For each subproblem, we calculate the corresponding variable by fixing the rest ones. In particular, for the hierarchical tensor ring decomposition problem, we update the variable from the first layer to the second layer. Experimental results  on color images and hyperspectral images (HSI) show that our method outperforms state-of-the-art algorithms in terms of peak signal-to-noise ratio (PSNR), structural similarity index (SSIM), relative square error (RSE) and spectral angle mapper (SAM).

The main contributions of this paper can be summarized as follows: 

1) We propose a hierarchical tensor-ring decomposition. The tensor ring factors are further decomposed by t-SVD in the new tensor network. This is the first tensor network using multiple kinds of operators as connections for factors. In fact, the t-SVD for factors in the second layer can be applied to other tensor networks, e.g., tensor train, and the general hierarchical tensor networks can be obtained.

2) We adopt the newly proposed tensor network based rank and total variation terms to simultaneously utilize the multidimensional global data structure and the locally piecewise smoothness data structure.

3) Experimental results on color images and HSI show that the proposed method outperforms state-of-the-art algorithms in terms of recovery accuracy.

The rest of this paper is organized as follows. In Section II, we give the notations and preliminaries about tensor ring decomposition. The proposed model for tensor completion is given in Section III. In Section IV, we present detailed solutions. In Section V, numerical results are demonstrated. The conclusion is drawn in Section VI.
\section{Notations and  Preliminaries}
\subsection{Notations}
In this paper, a scalar is denoted by standard lower case letter or uppercase letter, e.g., $a\in \mathbb{R},$ and a vector is denoted by boldface lowercase letter e.g., $\mathbf{a}\in \mathbb{R}^{I}$. A matrix is denoted by boldface capital letter, e.g., $\mathbf{A}\in \mathbb{R}^{I_{1}\times I_{2}}$. A tensor of order $N\geq 3$ is denoted by boldface calligraphic letters, e.g., $\mathbf{\mathcal{A}}\in \mathbb{R}^{I_{1}\times I_{2}\times\cdots\times I_{N}}$. The elements of tensor $\mathcal{A}\in \mathbb{R}^{I_{1}\times I_{2}\times\cdots\times I_{N}}$ is defined by $a_{i_{1},i_{2},i_{3},\cdots,i_{N}} $, where $i_{n}$ shows $n^{th}$ index of tensor $\mathcal{A}$. 

The Frobenius norm of $\mathcal{A}$ is defined by $\|\mathcal{A}\|_\text{F}=$ $\sqrt{\langle\mathcal{A},\mathcal{A}\rangle}$. The inner product of two tensors $\mathcal{A}, \mathcal{B}$ with the same size $\mathbb{R}^{I_{1}\times I_{2}\times\cdots\times I_{N}}$ can be defined as  $\langle\mathcal{A}, \mathcal{B}\rangle=\sum_{i_{1}} \sum_{i_{2}} \cdots \sum_{i_{N}} a_{{i_{1} i_{2} \ldots i_{N}}} b_{{i_{1} i_{2} \ldots i_{N}}}.$
 Letting $\mathbf{\mathcal{A}} \in$ $\mathbb{R}^{I_1 \times I_2 \times \cdots \times I_N}$ be an $N$-order tensor, the standard mode-$n$ unfolding of $\mathbf{\mathcal{A}}$ can be defined as  $\mathbf{A}_{(n)} \in \mathbb{R}^{I_{n} \times I_{1} \cdots I_{n-1} I_{n+1} \cdots I_{N}}.$
 Another mode-$n$ unfolding of tensor is often used in tensor ring, which is defined as $\mathbf{A}_{<n>} \in \mathbb{R}^{I_{n} \times I_{n+1} \cdots I_{N} I_{1} \cdots I_{n-1}} .$

The matrix nuclear norm of $\mathbf{A}$ is denoted as  $\|\mathbf{A}\|_{*}=\sum_{n=1}^{N} \sigma_{n}(\mathbf{A}),$ where $\sigma_{n}(\mathbf{A})$ is
the singular value of matrix $\mathbf{A}$. $d^{*}$ denotes the conjugate of $d.$ 

The 3D total variation of  $\mathbf{\mathcal{A}}\in \mathbb{R}^{I_1 \times I_2 \times I_3}$ can be formulated as follows:
$$\|\mathbf{\mathcal{A}}\|_{\mathrm{TV}}=\sum_{i_1=1}^{I_1}\sum_{i_2=1}^{I_2}\sum_{i_3=1}^{I_3}\left\|\nabla_{x}\mathbf{\mathcal{A}}\right\|_{1}+\left\|\nabla_{y}\mathbf{\mathcal{A}}\right\|_{1}+\left\|\nabla_{z}\mathbf{\mathcal{A}}\right\|_{1} $$
where $(\nabla_{x}\mathbf{\mathcal{A}})=\mathbf{\mathcal{A}}_{i_1,i_2,i_3}-\mathbf{\mathcal{A}}_{i_1+1,i_2,i_3}, (\nabla_{y}\mathbf{\mathcal{A}})=\mathbf{\mathcal{A}}_{i_1,i_2,i_3}-\mathbf{\mathcal{A}}_{i_1,i_2+1,i_3}, (\nabla_{z}\mathbf{\mathcal{A}})\mathbf{\mathcal{A}}_{i_1,i_2,i_3}=\mathbf{\mathcal{A}}_{i_1,i_2,i_3+1}$


\subsection{Preliminaries on tensor ring decomposition}

\begin{figure}[t]
	\includegraphics[scale=0.3]{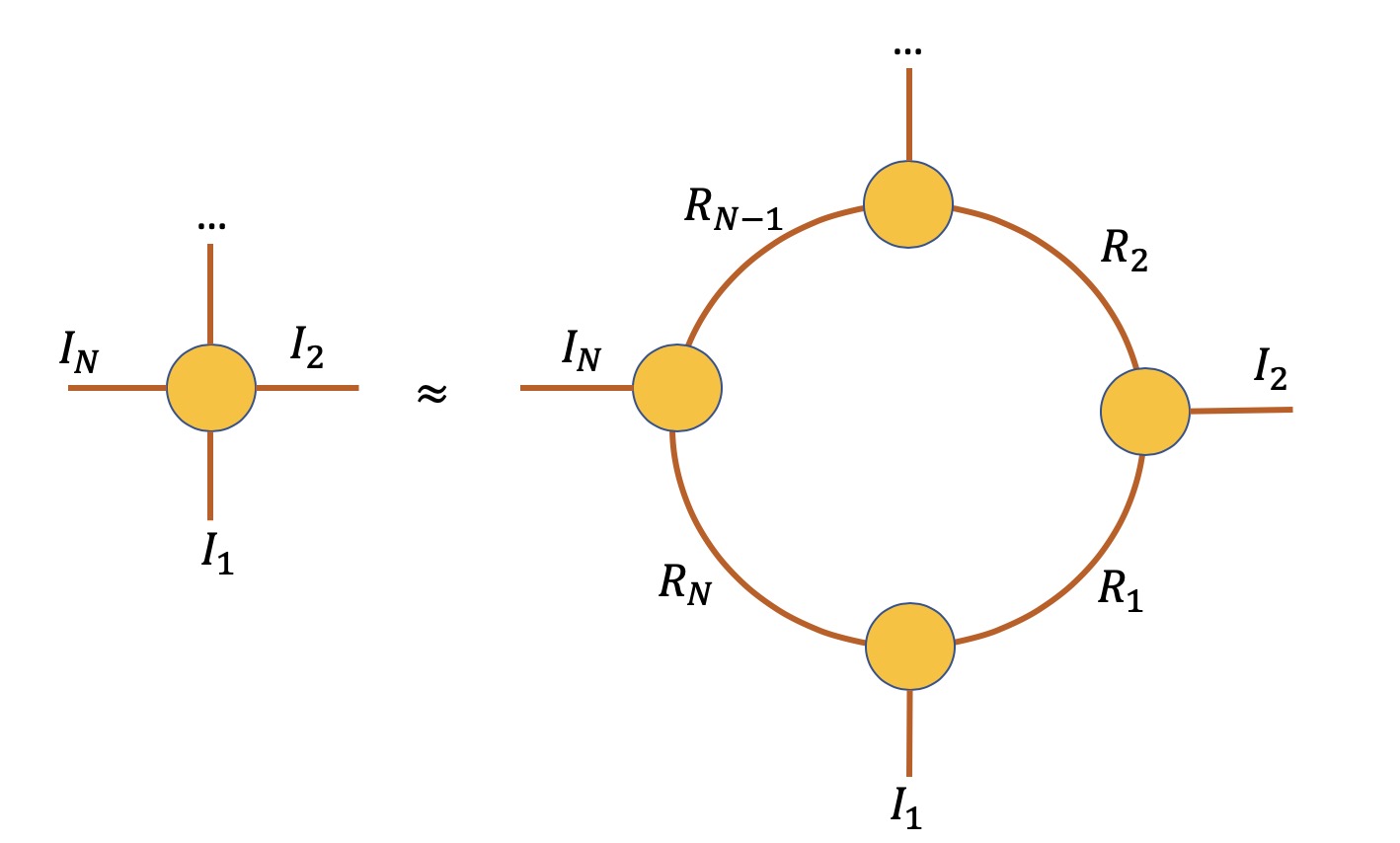}
	\caption{Graphical representation of tensor-ring (TR) model.}
	\label{16}
\end{figure}

\begin{definition}(tensor ring decomposition)\cite{zhao2016tensor}

Tensor ring decomposition represents a high-order tensor into multilinear products of low-order tensors in a circular form, where low-order tensors are called TR factors. The element-wise relationship of TR decomposition and generated tensors can be defined as follows: 
 \begin{equation}\label{d}
\mathcal{A}\left(i_{1}, i_{2}, \ldots, i_{N}\right)=\text { Trace }\left\{\prod_{n=1}^{N} \mathbf{\mathcal{G}}_{n}(i_n)\right\}
\end{equation}
where $ \text{Trace}\{\cdot\}$ is the matrix trace operation,  $\mathbf{\mathcal{G}}_{n}(i_n)\in$ $\mathbb{R}^{R_{n} \times R_{n+1}}$ is the  $i_{n}$$\text {-th } \text { mode-2 slice matrix of } \mathbf{\mathcal{G}}_{n},$ which can also be denoted by  $\mathbf{\mathcal{G}}_{n}\left(:, i_{n}, :\right)$ according to MatLab notation. Fig. \ref{16} gives a more intuitive representation of tensor ring decomposition.  For simplification, we use the notion $\mathcal{F}(\mathcal{G}_1,\cdots,\mathcal{G}_N)$ to represent the tensor ring decomposition of an N-order tensor.
\end{definition}

\begin{definition}(tensor ring rank)\cite{zhao2016tensor}

 In a tensor ring, the tensor ring rank is a vector $[R_1; R_2; \cdots; R_N]$. We set all TR-ranks to be equal in this paper, i.e. ${R_n}={R}$, $\forall n=1,\dots,N.$ 
\end{definition}

\begin{definition}(tensor permutation) \cite{wang2017efficient}

The $k^{th}$ tensor permutation of an tensor $\mathbf{\mathcal{A}} \in$ $\mathbb{R}^{I_1 \times I_2 \times \cdots \times I_N}$ can be defined as $\mathbf{\mathcal{A}}^{P_k}\in$ $\mathbb{R}^{{I_{1} \times I_{k+1} \times \cdots \times I_{N} \times I_{1} \times I_{2} \times \cdots \times I_{k-1}}} $ such that $\forall$ $k,i_k \in [1,I_k]$,
$$\mathbf{\mathcal{A}}^{P_k}\left(i_{k}, \cdots, i_{N}, i_{1}, \cdots, i_{k-1}\right)=\mathbf{\mathcal{A}}\left(i_{1}, \cdots, i_{N}\right)$$
and we can have the following result:
\end{definition}
\begin{lemma}
If $\mathbf{\mathcal{A}}=\mathcal{F}(\mathcal{G}_1,\cdots,\mathcal{G}_N),$  $\mathbf{\mathcal{A}}^{P_k}=$$\mathcal{F}\left(\mathcal{G}_{k} \mathcal{G}_{k+1} \cdots \mathcal{G}_{N} \mathcal{G}_{1} \cdots \mathcal{G}_{k-1}\right)$.
\end{lemma}

\begin{definition}
(tensor connect product)\cite{wang2017efficient}
 Assuming $\mathbf{\mathcal{G}}_{n}\in$ $\mathbb{R}^{R_{n-1}\times I_{n}\times R_{n}}$, the tensor connect product between $\mathbf{\mathcal{G}}_{n}$ and $\mathbf{\mathcal{G}}_{n+1}$ can be defined as:
$$\mathbf{\mathcal{G}}_{n} \mathbf{\mathcal{G}}_{n+1}\in\mathbb{R}^{R_{n-1}\times (I_{n}I_{n+1})\times R_{n}}$$
$$ =\text{reshape}\left(\mathbf{L}\left(\mathbf{\mathcal{G}}_{n}\right) \times \mathbf{R}\left(\mathbf{\mathcal{G}}_{n+1}\right)\right)$$
\newline Thus, tensor connect product can be formulated as follows:
$$\mathbf{\mathcal{G}}=\mathbf{\mathcal{G}}_{1}\cdots\mathbf{\mathcal{G}}_{N} \in \mathbb{R}^{R_{0}\times (I_{1}\cdots I_{N})\times R_{N}}.$$
\end{definition}
\begin{figure}[h]
	\centering
	\includegraphics[scale=0.6]{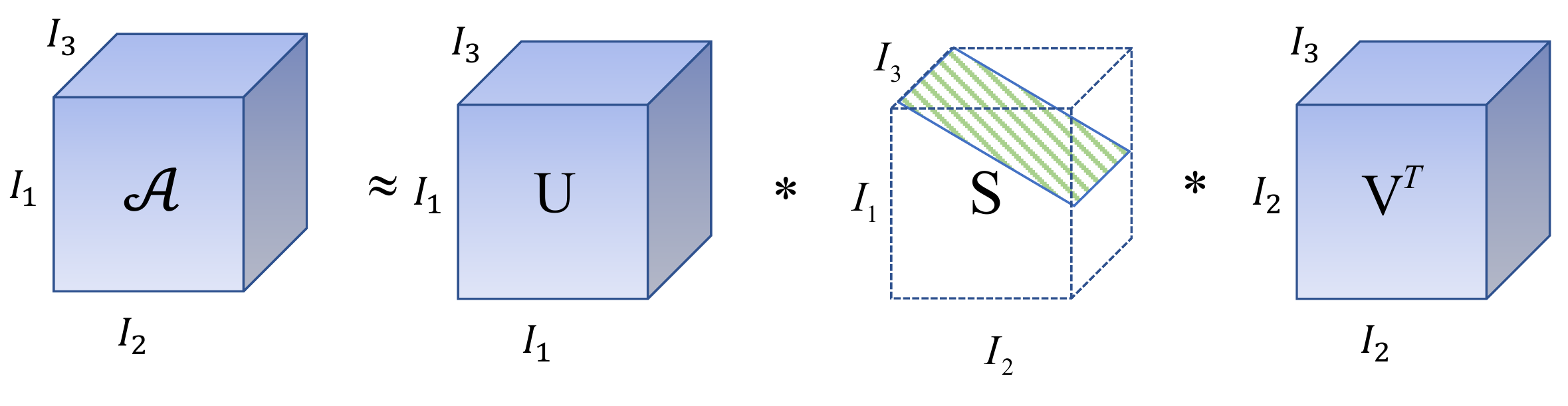}
	\caption{Illustration of t-SVD of  $\mathbf{\mathcal{A}}\in \mathbb{R}^{I_1\times I_2\times I_3}$}
	\label{tsvd}
\end{figure}
\subsection{Preliminaries on tensor singular value decomposition}
\begin{definition}(t-SVD)~\cite{zheng2018tensor}

The t-SVD of a $3\text{-order}$ tensor  $\mathbf{\mathcal{A}}\in \mathbb{R}^{I_1\times I_2\times I_3}$ can be represented as follows: 
\begin{equation}\label{SVD}
\mathbf{\mathcal{A}}=\mathbf{\mathcal{U}}\ast\mathbf{\mathcal{S}}\ast\mathbf{\mathcal{V}}^\text{T}
\end{equation}
Where $\mathbf{\mathcal{U}}\in \mathbb{R}^{I_1\times I_1\times I_3}$ and  $\mathbf{\mathcal{V}}\in \mathbb{R}^{I_2\times I_2\times I_3}$ are the orthogonal tensors, and  $\mathbf{\mathcal{S}}\in \mathbb{R}^{I_1\times I_2\times I_3}$ is an f-diagonal tensor. Fig. \ref{tsvd} illustrates the t-SVD of tensor $\mathbf{\mathcal{A}}$. 
\end{definition}
\begin{definition} (tensor tubal-rank)\cite{zheng2018tensor}
Let $\mathbf{\mathcal{A}}\in \mathbb{R}^{I_1\times I_2\times I_3}$ be a $3\text{-order}$ tensor. The tubal-rank of tensor $\mathbf{\mathcal{A}}$	denoted as $\operatorname{rank}_{\text{tubal}}(\mathbf{\mathcal{A}})$, is defined as the number of non-zero tubes of $\mathbf{\mathcal{S}}$. Where $\mathbf{\mathcal{S}}$ is defined in equation ( \ref{SVD}).
\end{definition}
\begin{definition}(tensor nuclear norm (TNN))\cite{zheng2018tensor}
	The nuclear norm of a tensor $\mathbf{\mathcal{A}}\in \mathbb{R}^{I_1\times I_2\times I_3}$ is defined as the sum of singular values of each frontal slice of $\bar{\mathbf{\mathcal{A}}}$ as follows:
	\begin{equation}
\|\mathbf{\mathcal{A}}\|_{\text{TNN}}\nonumber\\:=\sum_{i_3=1}^{I_3}\|\bar{\mathbf{A}}^{(i_3)}\|_{*}
	\end{equation}
where $\bar{\mathbf{A}}^{(i_3)}$ is the $n\text{-th}$ frontal slice of $\bar{\mathbf{\mathcal{A}}}$, and  $\bar{\mathbf{\mathcal{A}}}=\operatorname{fft}(\mathbf{\mathcal{A}}, [], 3)$.
\end{definition}

\section{Optimization Model}

\begin{figure}[t]
	\includegraphics[scale=0.3]{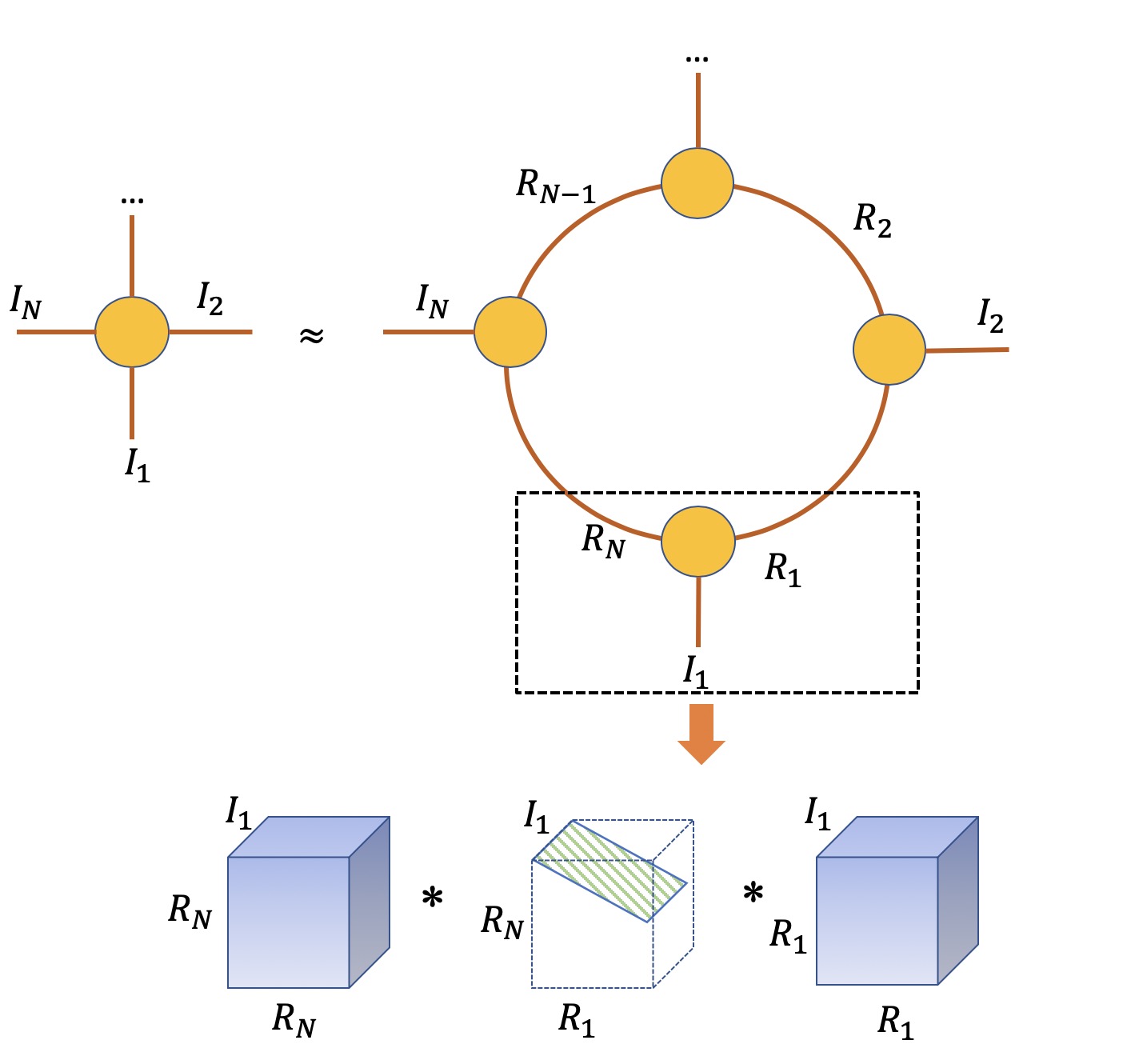}
	\caption{Graphical representation of a hierarchical tensor-ring (TR) model.}
	\label{fig:2}
\end{figure}

In this section, we propose a smooth low rank hierarchical tensor ring approximation (SHTRA)  for image completion. 
First, we develop a hierarchical tensor ring decomposition for more compact multi-way data representation.  
In this new decomposition, for the first layer,  traditional tensor ring decomposition is used to factorize a tensor into several 3-order tensors. For the second layer, each 3-order tensor can be further decomposed by t-SVD. Fig. \ref{fig:2} shows the details of hierarchical tensor-ring decomposition.

 The low rank approximation based the proposed hierarchical tensor-ring decomposition of $\mathcal{X}$ can be formulated as follows:
 \begin{equation}
 \min_{\mathcal{G}_{n,n=1,\cdots,N}} \frac{1}{2}\|\mathcal{X}-\mathcal{F}(\mathcal{G}_1, \cdots, \mathcal{G}_N)\|^{2}_{\operatorname{F}}+ \sum_{n=1}^{N}\text{rank}_{\text{tubal}}(\mathcal{G}_n)
 \end{equation}

To further enhance the recovery performance in many data processing, we can add TV term to exploit the piecewise smoothness structure. Thus, the optimization model of smooth low rank hierarchical  tensor-ring approximation for image completion can be written as:
\begin{eqnarray}\label{one}
&&\min_{\mathcal{X},~\mathcal{G}_{n, n=1,\cdots,N}}\frac{1}{2}\|\mathcal{X}-\mathcal{F}(\mathcal{G}_1, \cdots, \mathcal{G}_N)\|^{2}_{\operatorname{F}}+\lambda \text{TV}(\mathcal{X})\nonumber\\
&&+ \sum_{n=1}^{N}\text{rank}_{\text{tubal}}(\mathcal{G}_{n})\nonumber\\
&&\text{s. t.}\qquad \mathcal{P}_\mathbb{O}(\mathcal{X})=\mathcal{P}_\mathbb{O}(\mathcal{T})
\label{opt_model}
\end{eqnarray}
where $\mathcal{X}$ is the recovered low-rank tensor, $\mathcal{T}$ is the tensor for measurements, and  $\lambda$ is the trade-off parameter between hierarchical  tensor ring term and  TV term. $\mathcal{P}_\mathbb{O}(\mathcal{X})$ represents the observed entries. 

The tubal rank is non-convex, and one of its convex  surrogate is tensor nuclear norm. $\text{TV}\left(\mathbf{\mathcal{X}}\right)$ is a total variation of tensor $\mathbf{\mathcal{X}}$ which can be denoted as
 $\|\mathbf{\mathcal{D}\left(\mathbf{\mathcal{X}}\right)}\|_{1}$, i.e. the tensor total variation is the $ \ell_1 $ norm of all the differences along all the modes, where $ \mathcal{D}(\mathcal{X}) $ takes the differences.  In this way, the convex optimization model for (\ref{one}) can be formulated as
\begin{eqnarray}\label{second}
&&\min_{\mathcal{X},\mathcal{G}_{n, n=1,\cdots,N}}\frac{1}{2}\|\mathcal{X}-\mathcal{F}(\mathcal{G}_1, \cdots, \mathcal{G}_N)\|^{2}_{\operatorname{F}}+ \lambda \|\mathbf{\mathcal{D}}\left(\mathbf{\mathcal{X}}\right)\|_{1}\nonumber\\
&&+ \sum_{n=1}^{N}\|\mathbf{\mathcal{G}}_{n}\|_{\text{TNN}}\nonumber\\
&&\text{s. t.}\qquad \mathcal{P}_\mathbb{O}(\mathcal{X})=\mathcal{P}_\mathbb{O}(\mathcal{T}),
\end{eqnarray}
 where $\|.\|_{\text{TNN}}$ denotes the tensor nuclear norm. To solve the optimization model (\ref{second}), we introduce  additional tensor variables $\mathbf{\mathcal{Z}}$ and $\mathbf{\mathcal{M}}$ with the same size as $\mathbf{\mathcal{X}}$, and $\mathbf{Y}$ for tensor difference. The  equivalence of ((\ref{second}) can be obtained as follows:
\begin{equation} \label{a}
\min_{\mathbf{\mathcal{X}},\mathbf{\mathcal{G}}_{n}}\frac{1}{2}\|\mathbf{\mathcal{X}}-\mathcal{F}\left(\mathbf{\mathcal{G}}_1, \cdots, \mathbf{\mathcal{G}}_N\right)\|_\text{F}^2 + \lambda \|\mathbf{Y}\|_{1} +\sum_{n=1}^{N}\|\mathbf{\mathcal{G}}_{n}\|_{\text{TNN}}\nonumber\\ ,$$$$ \quad \text{s.t.}\quad \mathbf{\mathcal{P}}_\mathbb{O}\left(\mathbf{\mathcal{X}}\right)=\mathbf{\mathcal{P}}_\mathbb{O}\left(\mathbf{\mathcal{T}}\right), \mathbf{\mathcal{Z}}=\mathbf{\mathcal{X}}, \mathbf{Y}=\mathbf{\mathcal{D}}\left(\mathbf{\mathcal{Z}}\right), \sum_{n=1}^{N}\mathbf{\mathcal{M}}_{n}=\mathbf{\mathcal{G}}_{n} 
\end{equation}
Under the ADMM framework, this problem further can be converted into following form:
\begin{equation} \label{b}
\min_{\mathbf{\mathcal{X}},\mathbf{\mathcal{G}}_{n},\mathbf{\mathcal{Z}},\mathbf{\mathcal{M},\mathbf{Y}}}
\frac{1}{2}\|\mathbf{\mathcal{X}}-\mathcal{F}\left(\mathbf{\mathcal{G}}_1, \cdots, \mathbf{\mathcal{G}}_N\right)\|_\text{F}^2 + \lambda \|\mathbf{Y}\|_{1} +\sum_{n=1}^{N}\|\mathbf{\mathcal{M}}_{n}\|_{\text{TNN}}\nonumber\\ $$$$+\left\langle\Lambda_{1},\mathbf{\mathcal{Z}}-\mathbf{\mathcal{X}}\right\rangle+{\frac{\beta_{1}}{2}\|\mathbf{\mathcal{Z}}-\mathbf{\mathcal{X}}\|_{\mathrm{F}}^{2}}\ +\left\langle\Lambda_{2},\mathbf{Y}-\mathbf{\mathcal{D}}\left(\mathbf{\mathcal{Z}}\right)\right\rangle+$$$${\frac{\beta_{2}}{2}\|\mathbf{Y}-\mathbf{\mathcal{D}}\left(\mathbf{\mathcal{Z}}\right)\|_{\mathrm{F}}^{2}}
+\sum_{n=1}^{N}\left(\left\langle\Lambda_{3}^{(n)},\mathbf{\mathcal{M}}_{n}-\mathbf{\mathcal{G}}_{n}\right\rangle+{\frac{\beta_{3}}{2}\|\mathbf{\mathcal{M}}_{n}-\mathbf{\mathcal{G}}_{n}\|_{\mathrm{F}}^{2}} \right)$$$$ \quad \text{s.t.}\quad \mathbf{\mathcal{P}}_\mathbb{O}\left(\mathbf{\mathcal{X}}\right)=\mathbf{\mathcal{P}}_\mathbb{O}\left(\mathbf{\mathcal{T}}\right)
\end{equation}
where $\beta_{1}, \beta_{2}, \beta_{3}$ are positive penalty scalars, and  $\Lambda_{1}, \Lambda_{2}, \Lambda_{3}$ are the dual variables. This problem can be solved by updating each variable with others fixed. 

The first sub-problem optimizes the variable $\mathbf{\mathcal{G}}$ with other variables fixed, which can be written as:
\begin{equation}\label{c}
\min_{\mathbf{\mathcal{G}}_{n}} \frac{1}{2}\|\mathbf{\mathcal{X}}-\mathcal{F}\left(\mathbf{\mathcal{G}}_{1},\cdots, \mathbf{\mathcal{G}}_N\right)\|_\text{F}^2\\$$$$ +\sum_{n=1}^{N}\left(\left\langle\Lambda_{3}^{(n)},\mathbf{\mathcal{M}}_{n}-\mathbf{\mathcal{G}}_{n}\right\rangle+{\frac{\beta_{3}}{2}\|\mathbf{\mathcal{M}}_{n}-\mathbf{\mathcal{G}}_{n}\|_{\mathrm{F}}^{2}}\right)
\end{equation}
The second sub-problem optimizes the variable $\mathbf{\mathcal{M}}$ and keeps  other variables fixed, can be written as:
\begin{equation}\label{g}
\min_{\mathbf{\mathcal{M}}} \sum_{n=1}^{N}\|\mathbf{\mathcal{M}}_{n}\|_{\text{TNN}}\\ +\sum_{n=1}^{N}\left\langle\Lambda_{3}^{(n)},\mathbf{\mathcal{M}}_{n}-\mathbf{\mathcal{G}}_{n}\right\rangle$$$$+ \sum_{n=1}^{N}{\frac{\beta_{3}}{2}\|\mathbf{\mathcal{M}}_{n}-\mathbf{\mathcal{G}}_{n}\|_{\mathrm{F}}^{2}}
\end{equation}
The third sub-problem optimizes the variable $\mathbf{\mathcal{Z}}$ while keeping the others fixed, can be written as:
\begin{equation}\label{k}
\min_{\mathbf{\mathcal{Z}}} \left\langle\Lambda_{1},\mathbf{\mathcal{Z}}-\mathbf{\mathcal{X}}\right\rangle+{\frac{\beta_{1}}{2}\|\mathbf{\mathcal{Z}}-\mathbf{\mathcal{X}}\|_{\mathrm{F}}^{2}}+\left\langle\Lambda_{2},\mathbf{Y}-\mathbf{\mathcal{D}}\left(\mathbf{\mathcal{Z}}\right)\right\rangle $$$$+{\frac{\beta_{2}}{2}\|\mathbf{Y}-\mathbf{\mathcal{D}}\left(\mathbf{\mathcal{Z}}\right)\|_{\mathrm{F}}^{2}}
\end{equation}
The fourth sub-problem on $\mathbf{Y}$ can be written as:
\begin{equation}\label{m}
\min_{\mathbf{Y}} \lambda \|\mathbf{Y}\|_{1}+ \left\langle\Lambda_{2},\mathbf{Y}-\mathbf{\mathcal{D}}\left(\mathbf{\mathcal{Z}}\right)\right\rangle+{\frac{\beta_{2}}{2}\|\mathbf{Y}-\mathbf{\mathcal{D}}\left(\mathbf{\mathcal{Z}}\right)\|_{\mathrm{F}}^{2}}
\end{equation}
The fifth sub-problem on $\mathcal{X}$ can be written as:
\begin{equation}\label{o}
\min_{\mathbf{\mathcal{X}}} \frac{1}{2}\|\mathbf{\mathcal{X}}-\mathcal{F}\left(\mathbf{\mathcal{G}}_1, \cdots, \mathbf{\mathcal{G}}_N\right)\|_\text{F}^2 + \left\langle\Lambda_{1},\mathbf{\mathcal{Z}}-\mathbf{\mathcal{X}}\right\rangle$$$$+{\frac{\beta_{1}}{2}\|\mathbf{\mathcal{Z}}-\mathbf{\mathcal{X}}\|_{\mathrm{F}}^{2}}\quad \text{s.t.}\quad \mathbf{\mathcal{P}}_\mathbb{O}\left(\mathbf{\mathcal{X}}\right)=\mathbf{\mathcal{P}}_\mathbb{O}\left(\mathbf{\mathcal{T}}\right)
\end{equation}
The sub-problem on updating dual variables can be written as: 
\begin{equation}\label{p}
\min_{\Lambda_{1},\Lambda_{2},\Lambda_{3}^{(n)},n=1,\cdots,N} \left\langle\Lambda_{1},\mathbf{\mathcal{Z}}-\mathbf{\mathcal{X}}\right\rangle+{\frac{\beta_{1}}{2}\|\mathbf{\mathcal{Z}}-\mathbf{\mathcal{X}}\|_{\mathrm{F}}^{2}}\ +\left\langle\Lambda_{2},\mathbf{Y}-\mathbf{\mathcal{D}}\left(\mathbf{\mathcal{Z}}\right)\right\rangle$$$$+{\frac{\beta_{2}}{2}\|\mathbf{Y}-\mathbf{\mathcal{D}}\left(\mathbf{\mathcal{Z}}\right)\|_{\mathrm{F}}^{2}}
+\sum_{n=1}^{N}\left\langle\Lambda_{3}^{(n)},\mathbf{\mathcal{M}}_{n}-\mathbf{\mathcal{G}}_{n}\right\rangle$$$$+\sum_{n=1}^{N}{\frac{\beta_{3}}{2}\|\mathbf{\mathcal{M}}_{n}-\mathbf{\mathcal{G}}_{n}\|_{\mathrm{F}}^{2}}
\end{equation}

\section{Solution}
\subsection{The solution of subrpoblem (\ref{c})}
The sub-problem (\ref{c}) can be reformulated as:
\begin{equation}\label{e}
\min_{\mathbf{\mathcal{G}}_{n}} \frac{1}{2}\|\mathbf{\mathcal{X}}_{<n>}-(\mathcal{G}_n)_{<2>}\left(\mathbf{\mathcal{G}}_{<2>}^{(\neq n)}\right)^\text{T}\|_\text{F}^2\\$$$$ +\sum_{n=1}^{N}\left(\left\langle\Lambda_{3}^{(n)},\mathbf{\mathcal{M}}_{n}-\mathbf{\mathcal{G}}_{n}\right\rangle+{\frac{\beta_{3}}{2}\|\mathbf{\mathcal{M}}_{n}-\mathbf{\mathcal{G}}_{n}\|_{\mathrm{F}}^{2}}\right)
\end{equation}
where  $\mathbf{\mathcal{G}}_{<2>}^{(\neq k)} \in \mathbb{R}^{\prod_{n=1, n\neq k}^{N}I_{i}\times R_{k-1}R_{k}}$ is a sub-chain dimension-mode unfolded matrix generated by merging all the core unfolded matrices except the $k\text{-th}$ core.   $(\mathcal{G}_n)_{<2>} \in \mathbb{R}^{I_{k}\times R_{k}R_{k-1}}$ is a $k \text{-th}$ core dimension-mode unfolded matrix. 

The solution for (\ref{e}) is as follows:
\begin{equation}\label{f}
	(\mathcal{G}_n)_{<2>}=\left( \mathbf{\mathcal{X}}_{<n>}\mathbf{\mathcal{G}}_{<2>}^{(\neq n)}+{\Lambda_{3}}_{<2>}^{(n)}+\beta_{3}(\mathcal{M}_n)_{<2>}\right)$$$$ \left(\left(\mathbf{\mathcal{G}}_{<2>}^{(\neq n)}\right)^{T} (\mathcal{G}_n)_{<2>}^{(\neq n)}+ \beta_{3} \mathbf{I} \right)^{-1}
\end{equation}
where $\mathbf{I}$ is an identity matrix.

\subsection{The solution of sub-problem (\ref{g})}
The optimization problem (\ref{g}) can be equivalent to: 
\begin{equation}\label{i}
\sum_{n=1}^{N}\left(\frac{1}{\beta_{3}}\|\mathbf{\mathcal{M}}_{n}\|_{\text{TNN}}\\+\frac{1}{2}\|\mathbf{\mathcal{M}}_{n}-\left(\mathbf{\mathcal{G}}_{n}-\frac{\Lambda_{3}^{(n)}}{\beta_{3}}\right)\|_\text{F}^{2}\right)
\end{equation}
Let $\tau=\frac{1}{\beta_{3}}$, $\mathbf{\mathcal{L}}_{n}={\mathcal{G}_{n}}-\frac{\Lambda_{3}^{(n)}}{\beta_{3}}$ and $\mathcal{M}_n \in \mathbb{R}^{I_1\times I_2\times I_3}$ be a $3^{rd}\text{-order}$ core tensor. Therefore, (\ref{i}) can be rewritten as:
\begin{equation}\label{j}
\sum_{n=1}^{N}\left(\tau \|\mathbf{\mathcal{M}}_{n}\|_{\text{TNN}}\\+\frac{1}{2}\|\mathbf{\mathcal{M}}_{n}-\mathbf{\mathcal{L}}_{n}\|_\text{F}^{2}\right)
\end{equation}
 (\ref{j}) can be computed by using tensor singular value thresholding (t-SVT) as follows \cite{lu1804tensor}:
\begin{equation}\label{M}
\mathcal{M}_{n}=\operatorname{SVT}_{\tau}(\mathcal{L}_n)
\end{equation}
where $n=1,\cdots,N$ and for each $\tau>0,$ tensor singular value thresholding (t-SVT) operator can be defined as:
\begin{equation}\label{DT}
\mathbf{\operatorname{SVT}}_{\tau}\left(\mathbf{\mathcal{L}}\right)=\mathbf{\mathcal{U}}\ast\mathbf{\mathcal{S}}_{\tau}\ast\mathbf{\mathcal{V}}^{\ast,}
\end{equation}
where 
\begin{equation}
\mathbf{\mathcal{S}}_{\tau}=\operatorname{ifft}\left(\left(\bar{\mathbf{\mathcal{S}}}-\tau\right)_{+},\left[ \right],3\right)
\end{equation}
where $\bar{\mathbf{\mathcal{S}}}$ is a real tensor . $t_+$ indicates the positive part i.e. $t_+=\max \left(t,0\right)$.  This operator applies the soft-thresholding rule to singular values $\bar{\mathbf{\mathcal{S}}}$ of each frontal slice of $\bar{\mathbf{\mathcal{L}}}$.
The detailed solutions are concluded in Algorithm \ref{algorithm:tsvt}.

\begin{algorithm}[t]
	\caption{Tensor Singular Value Thresholding (t-SVT)}
	\label{algorithm:tsvt}
	\begin{algorithmic}
		\STATE \textbf{Input:} $\mathcal{L}_n\in \mathbb{R}^{R_{n-1} \times R_n \times I_n}$, $\tau>0$
		\STATE  \textbf{Output:} $\mathbf{\operatorname{SVT}}_{\tau}\left(\mathbf{\mathcal{L}}\right)$ as defined in \ref{DT}
		\STATE Compute $\bar{\mathbf{\mathcal{L}}}=\operatorname{fft}\left(\mathbf{\mathcal{L}}, [], 3\right)$.
		\STATE Perform Matrix SVT on each frontal slice of $\bar{\mathbf{\mathcal{L}}}$ by 
		\FOR {$k=1, \cdots, \frac{I_n+ 1}{2}$}
		\STATE	$[U, S, V]=\operatorname{SVD}\left(\bar{\mathbf{L}}^{(k)}\right)$
		\STATE $\bar{\mathbf{W}}=U\left(S-\tau\right)_{+} V^{*};$
		\ENDFOR
		\FOR {$k=\left[\frac{I_n + 1}{2}\right]+1, \cdots, I_n$}
		\STATE $\bar{\mathbf{W}}^{(k)}=\operatorname{conj}\left(\bar{\mathbf{W}}^{(I_n-k + 2)}\right)$
		\ENDFOR
		\STATE Compute $\mathbf{\operatorname{SVT}}_{\tau}\left(\mathbf{\mathcal{L}}\right)=\operatorname{ifft}\left(\bar{\mathbf{W}}, [], 3\right)$
	\end{algorithmic} 
\end{algorithm}
\begin{figure}[htbp]
\begin{center}
\includegraphics[scale=0.30]{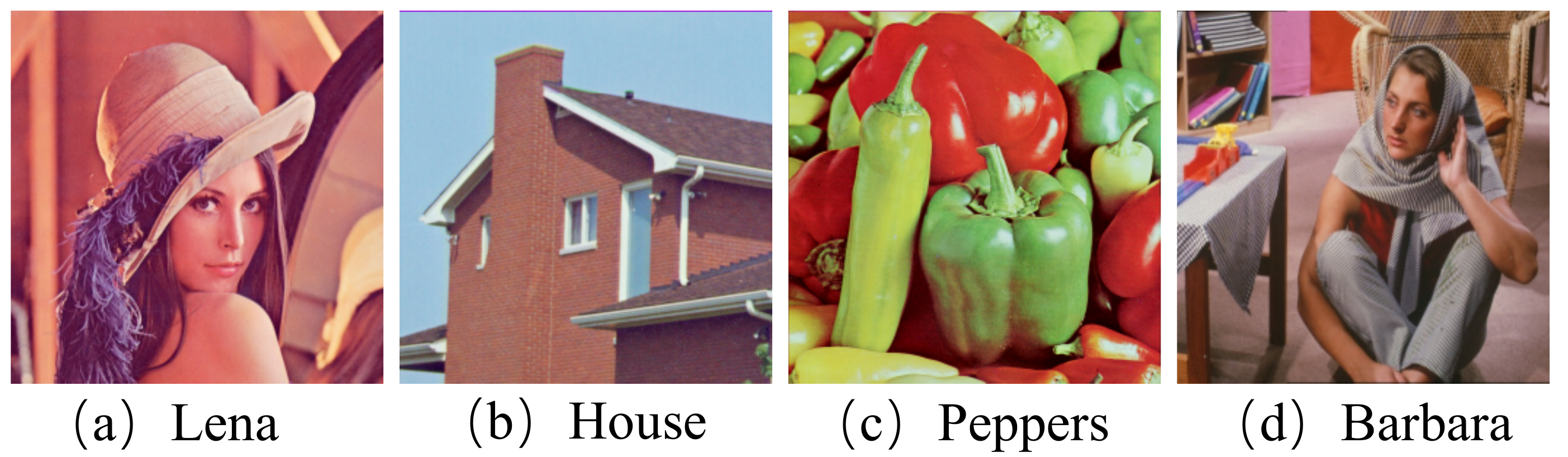}
\caption{Testing color images with the size $256 \times 256 \times 3$.}
\label{fig:A}
\end{center}
\end{figure}
\begin{figure}[htbp]
\begin{center}
\includegraphics[scale=0.30]{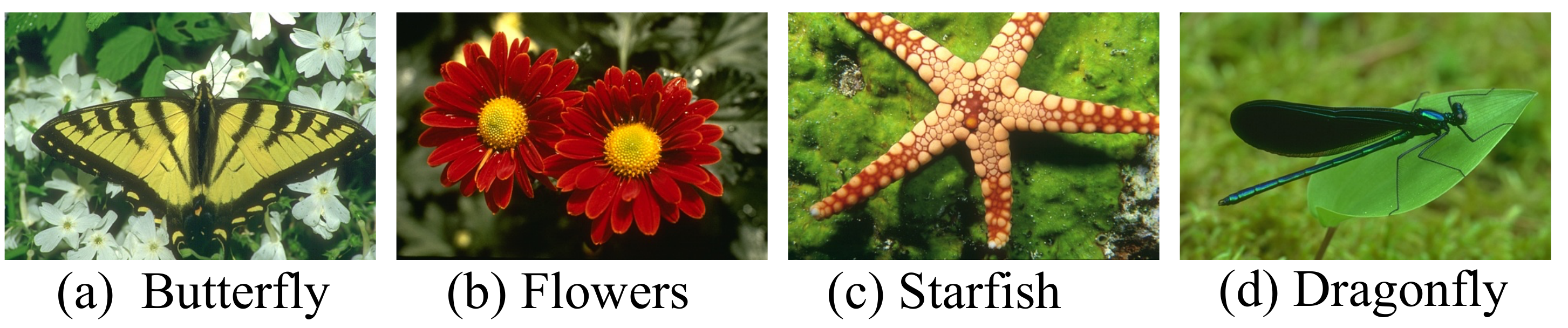}
\caption{Testing color images with the size $321 \times 481 \times 3$.}
\label{fig:B}
\end{center}
\end{figure}
\begin{figure}[htbp]
	\centering
	\includegraphics[scale=0.5]{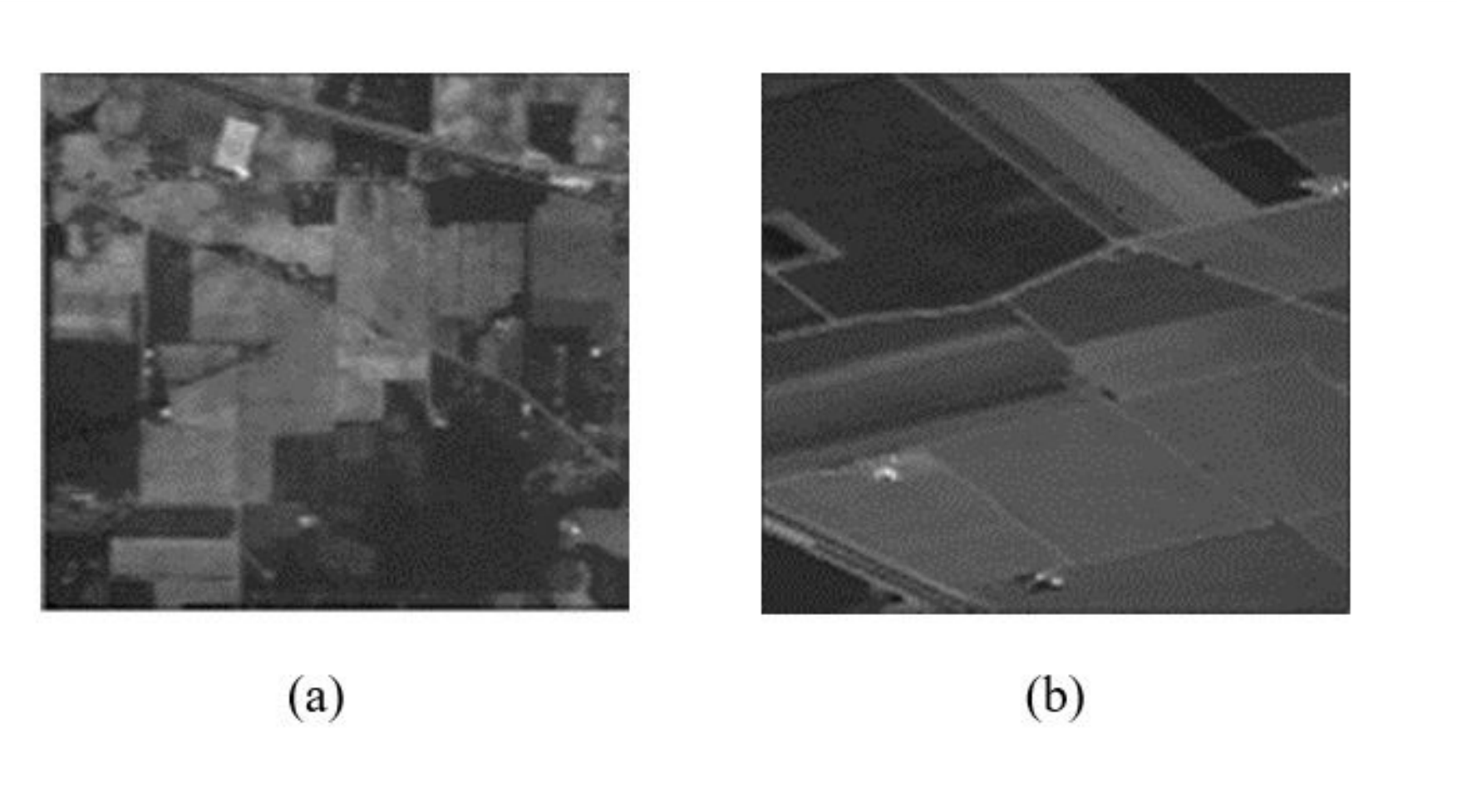}
	\caption{Hyperspectral images  ($10^{th}$ band grayscale) (a) Hyperspectral image 'Indianpines' with the size $145 \times 145 \times 224$. (b) Hyperspectral image 'Salinas' with size  $512 \times 217 \times 224$.}
	\label{2}
\end{figure}
\begin{figure*}[htbp]
	\centering
	\subfigure[PSNR vs SR]{
	\includegraphics[scale=0.36]{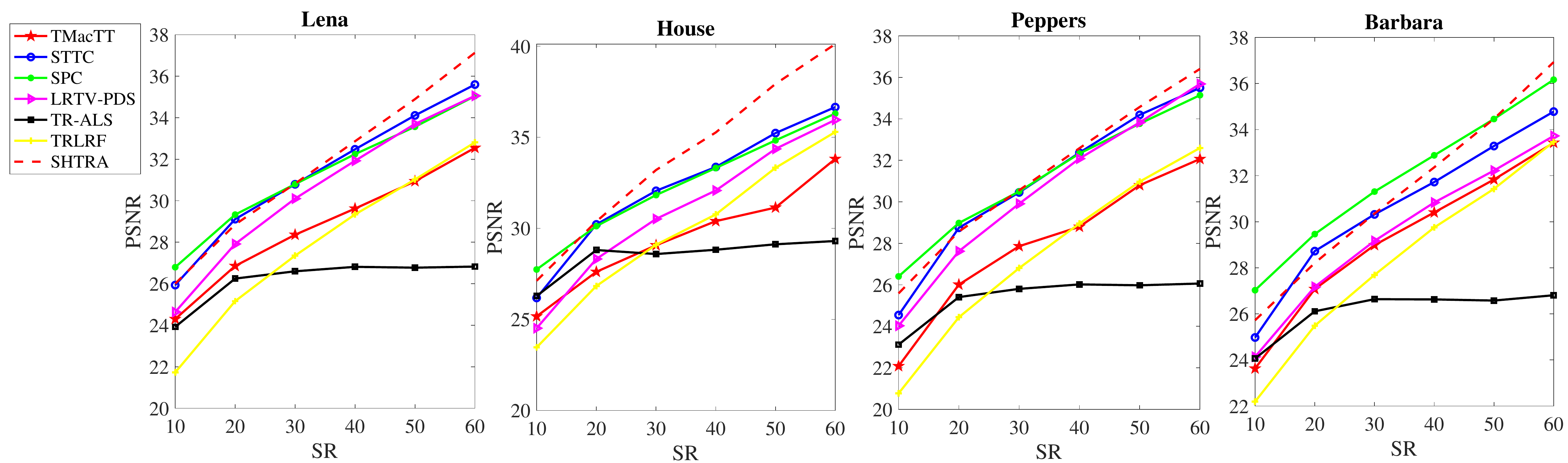}}\\
	\subfigure[SSIM vs SR]{
	\includegraphics[scale=0.36]{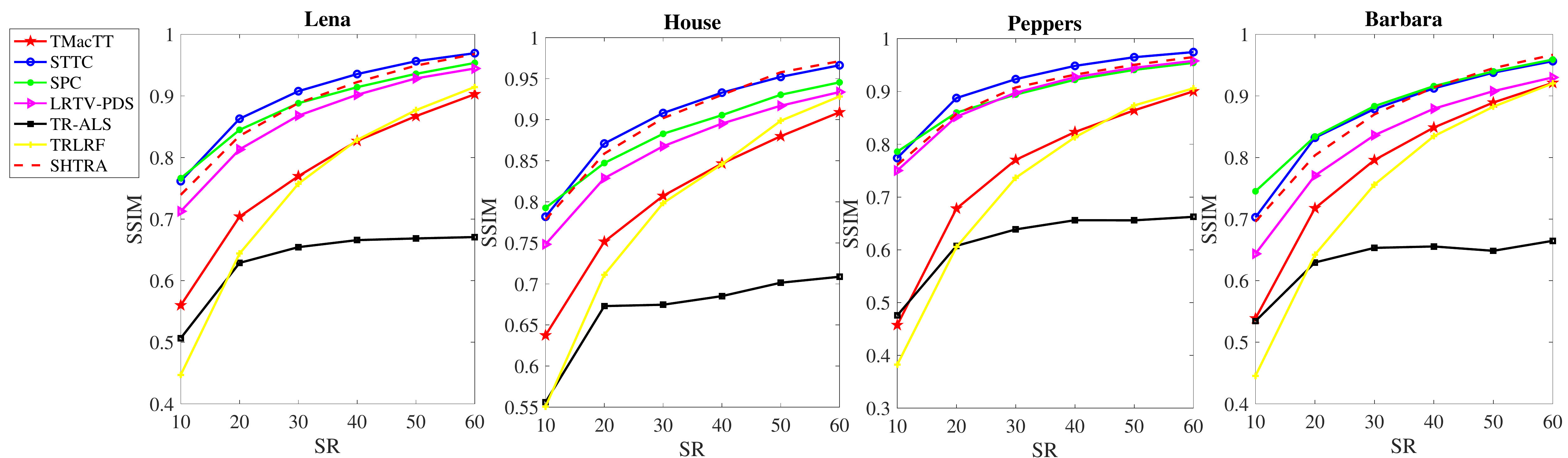}}	\\
	\subfigure[RSE vs SR]{
	\includegraphics[scale=0.36]{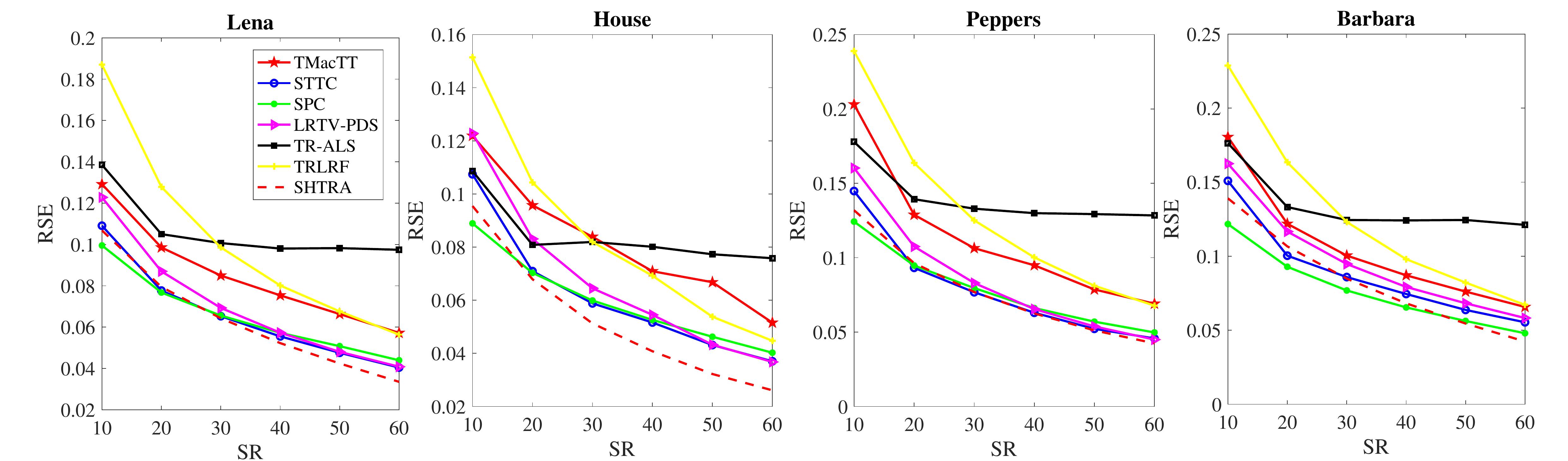}}	
	\caption{Comparison of different methods on color images with the size $256 \times 256 \times 3$ when different sampling ratios change from 10\% to 60\% in terns of PSNR, SSIM and RSE.}
	\label{fig:psnrA}
\end{figure*}

\begin{figure*}[htbp]
	\centering
	\includegraphics[scale=0.4]{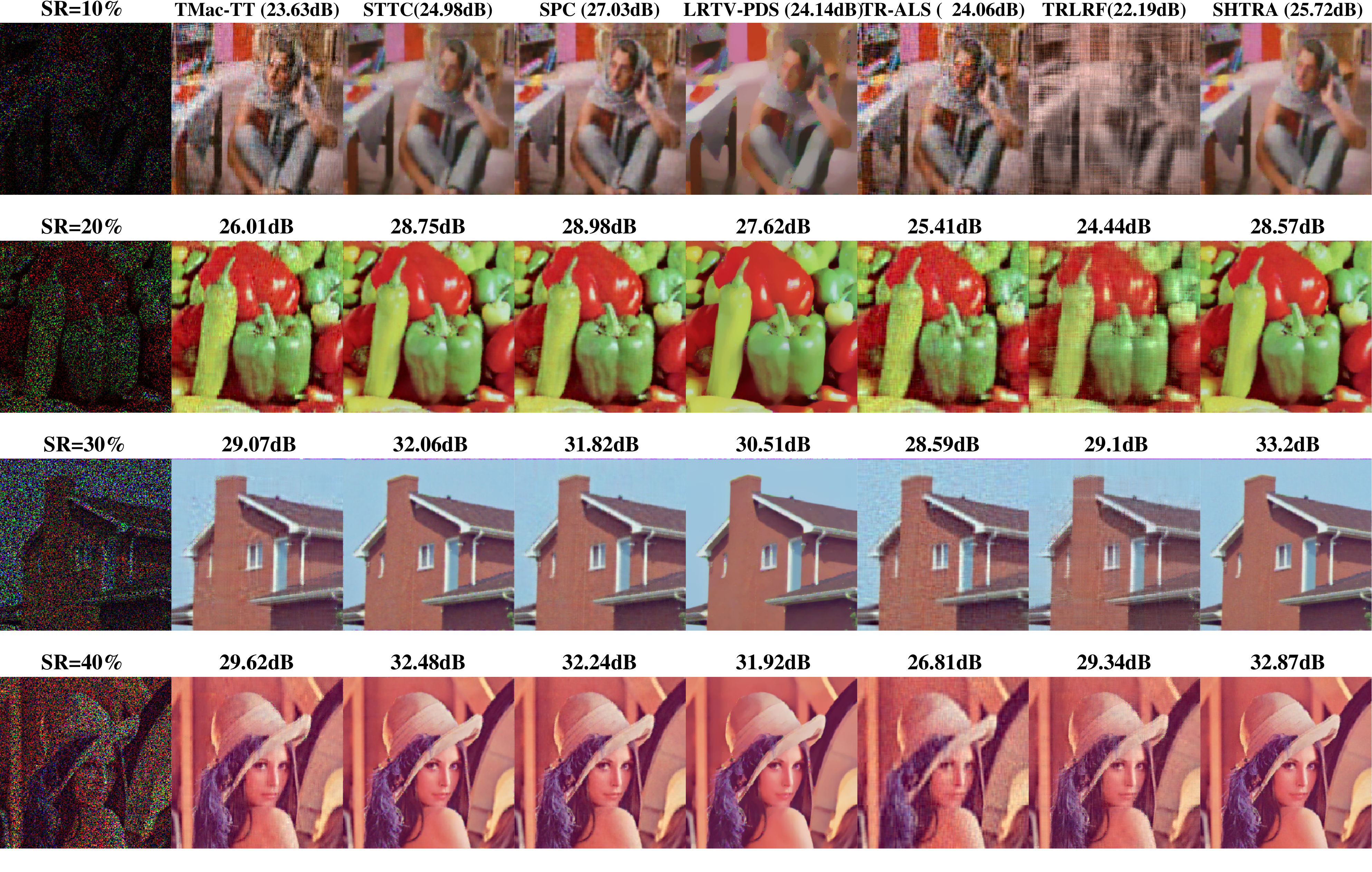}
	\caption{Examples on image completion using different methods with different sampling ratios. }
	\label{imageA}
\end{figure*}

\begin{figure*}[htbp]
	\centering
	\subfigure[PSNR vs SR]{
	\includegraphics[scale=0.36]{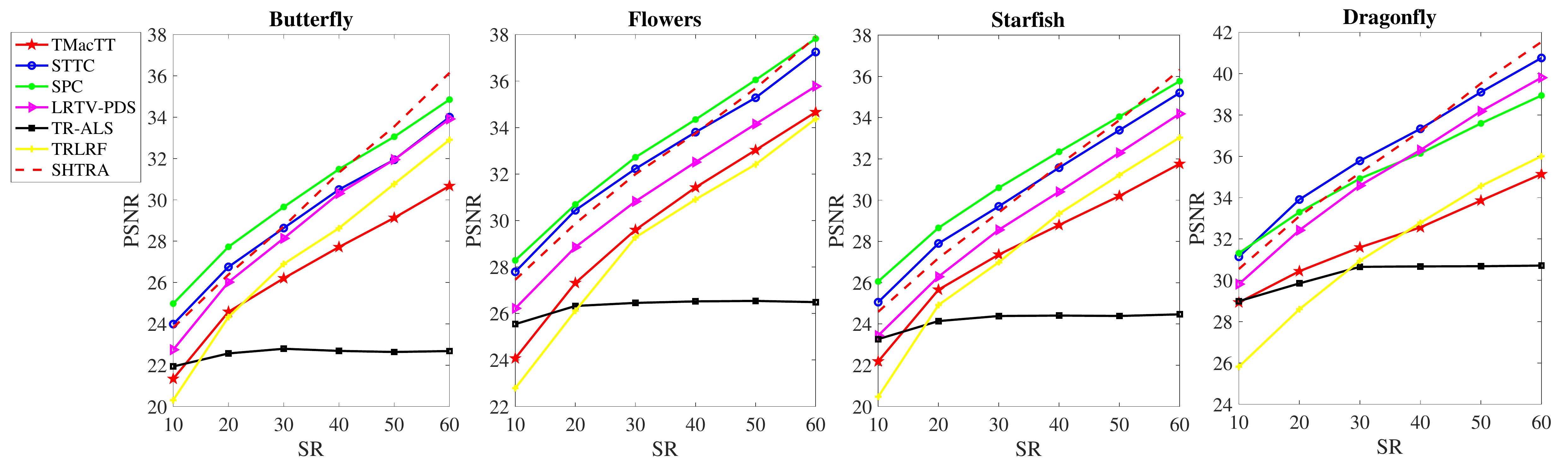}}\\
	\subfigure[SSIM vs SR]{
	\includegraphics[scale=0.36]{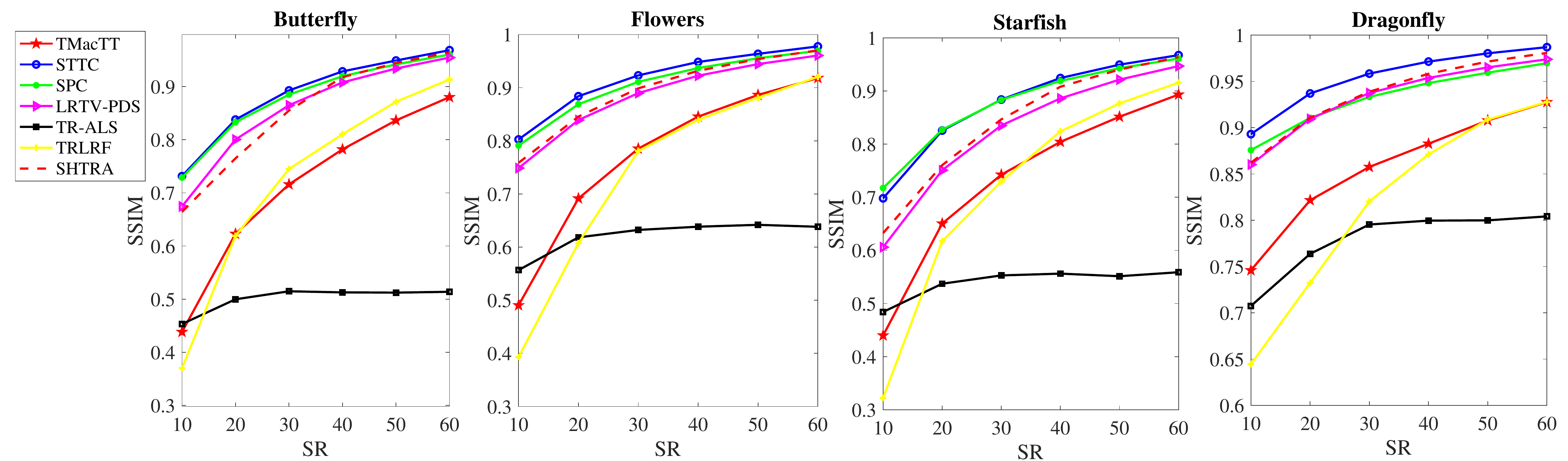}}	\\
	\subfigure[RSE vs SR]{
	\includegraphics[scale=0.36]{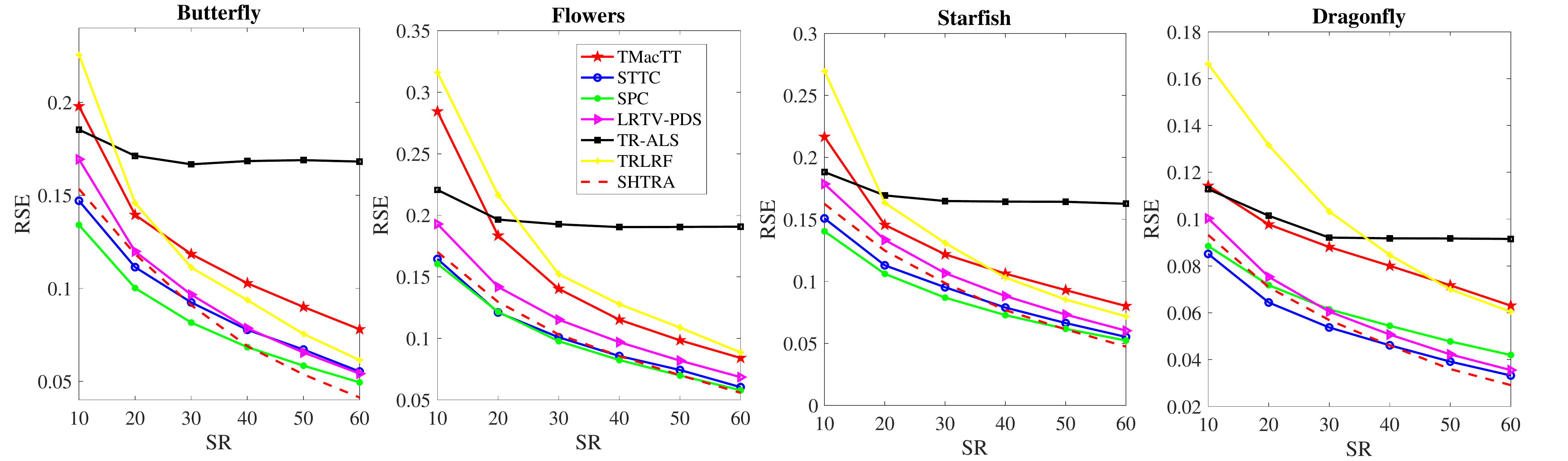}}	
	\caption{Comparison of different methods on color images with the size  $321 \times 481 \times 3$ when different sampling ratios change from 10\% to 60\% in terns of PSNR, SSIM and RSE.}
	\label{fig:psnrB}
\end{figure*}
\begin{figure*}[h]
	\centering
	\includegraphics[scale=0.4]{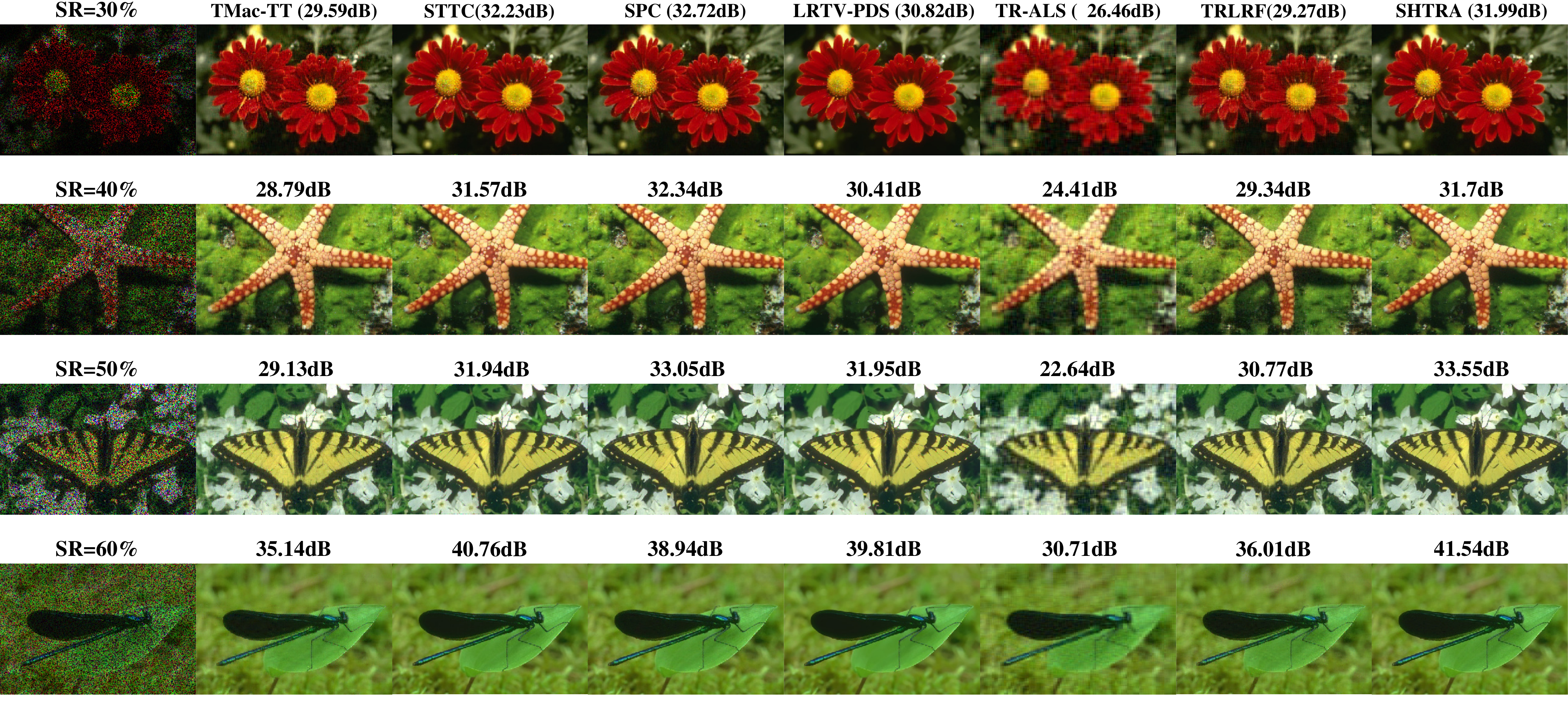}
	\caption{Examples on image completion using different methods with different sampling ratios.  }
	\label{imageB}
\end{figure*}
\begin{figure*}[h]
	\centering
	\includegraphics[scale=0.4]{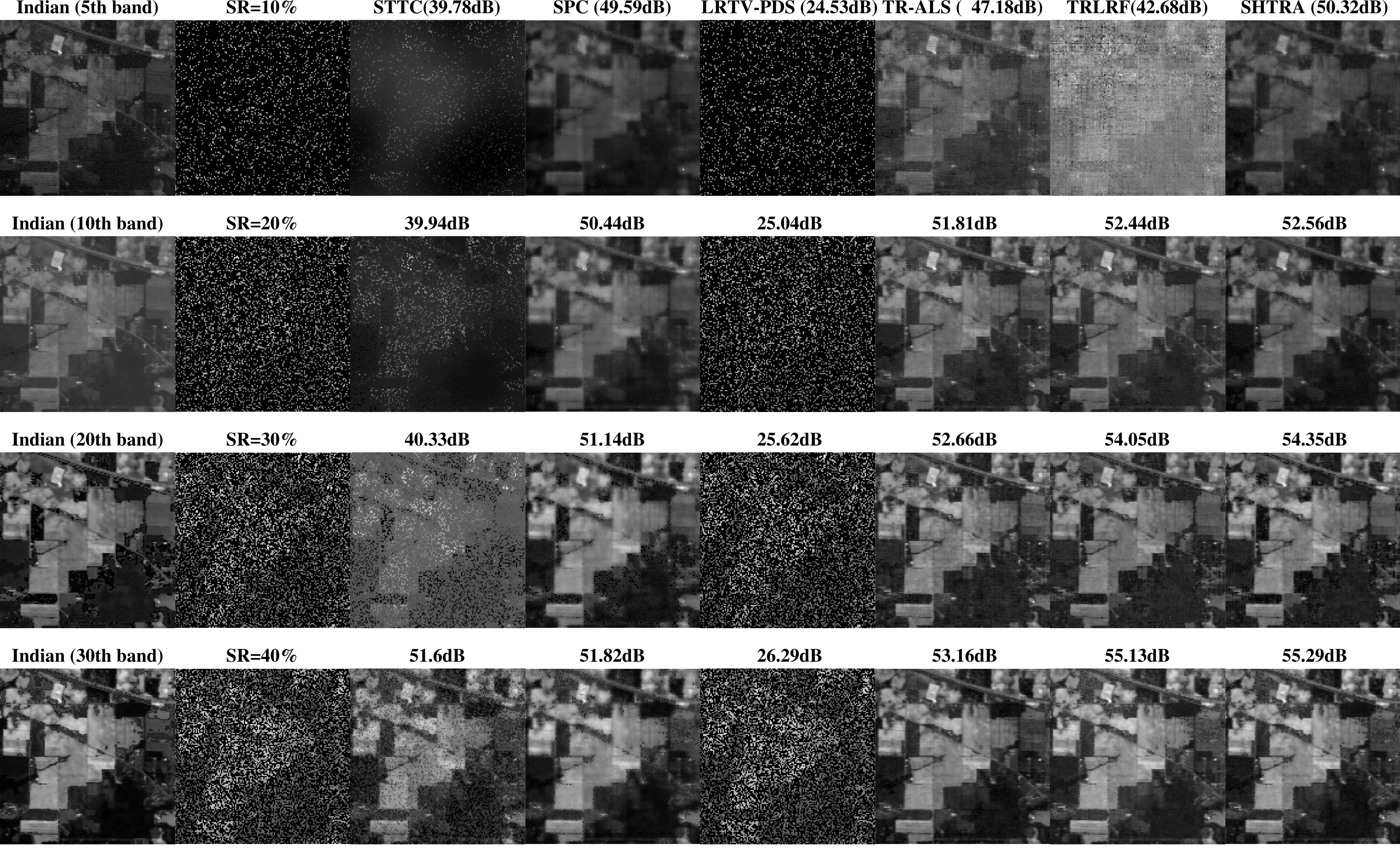}
	\caption{Examples on Indianpines completion using different methods with different sampling ratios. }
	\label{HSIIndian}
\end{figure*}
\begin{figure*}[h]
	\centering
	\includegraphics[scale=0.65]{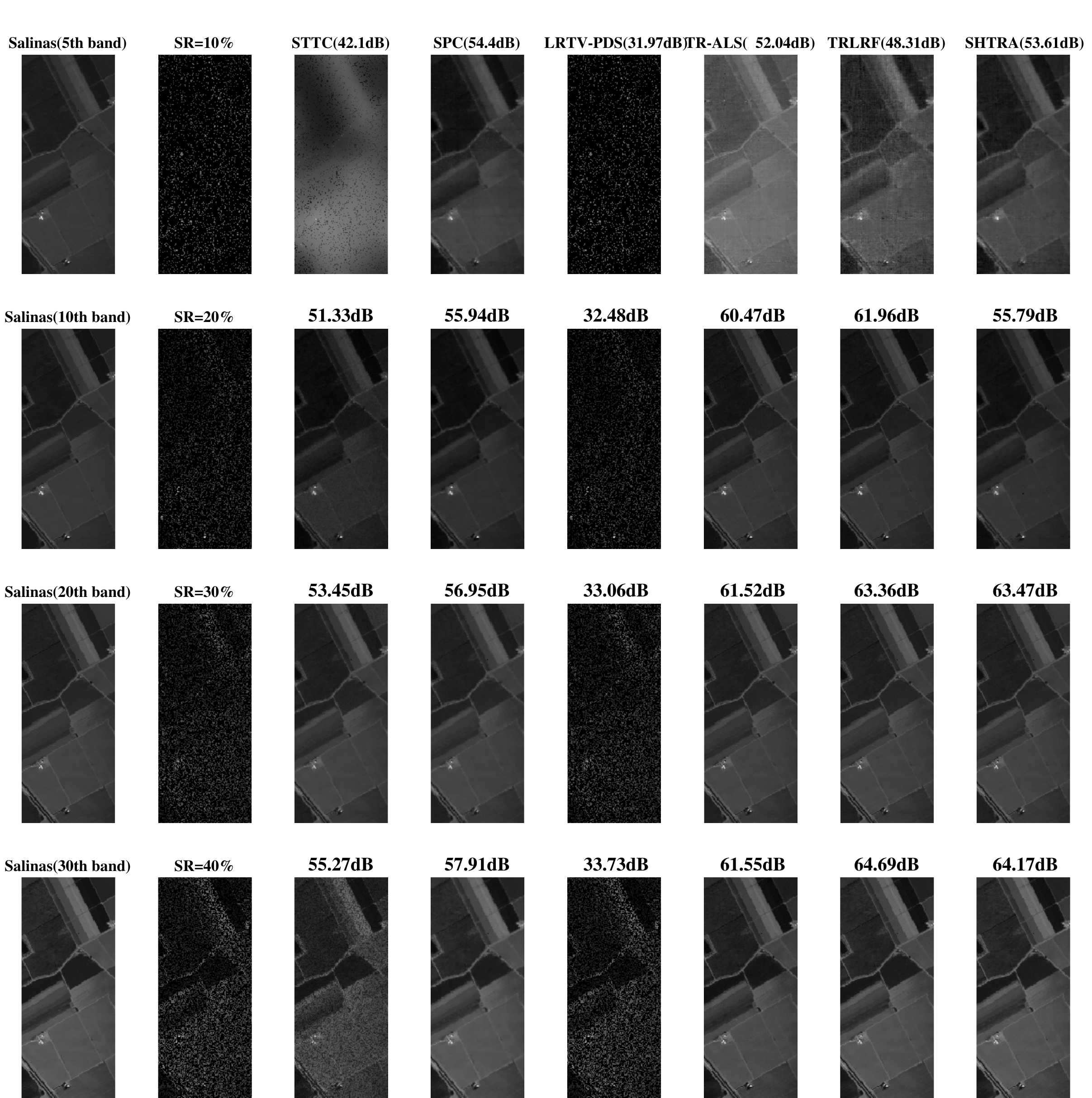}
	\caption{Examples on Salinas completion using different methods with different sampling ratios. }
	\label{HSISaline.eps}
\end{figure*}

\subsection{The solution of sub-problem (\ref{k})}
This sub-problem can be solved by differentiable function with respect to $\mathbf{\mathcal{Z}}$. The minimization condition is equivalent to the following linear equation \cite{8421084}:
\begin{equation}\label{l}
\mathbf{\mathcal{Z}}=\frac{\mathbf{\mathcal{D}}^{\ast}\left(\Lambda_{2}-\beta_{2}\mathbf{Y}\right)+\beta_{1}\mathbf{\mathcal{X}}-\Lambda_{1}}{\beta_{1}\mathbf{I}+\beta_{2}\mathbf{\mathcal{D}}^{\ast}\mathbf{\mathcal{D}}}
\end{equation}
where $\mathbf{\mathcal{D}}^{\ast}$ is adjoint of $\mathbf{\mathcal{D}}$. As the structure of operator  $\mathbf{\mathcal{D}}^{\ast}\mathbf{\mathcal{D}}$ is block circulant, it can be transformed into the Fourier domain and fast calculated. By using the off-the-shelf conjugate gradient technique, and the fast computation of $\mathbf{\mathcal{Z}}$ can be written as:
\begin{equation}\label{Z}
\mathbf{\mathcal{Z}}=\operatorname{ifftn}\left(\frac{\operatorname{fftn}\left(\mathbf{\mathcal{J}}\right)}{\beta_{1}\mathbf{I}+\beta_{2}\left(\operatorname{fftn}\left(\mathbf{\mathcal{D}}^{\ast}\mathbf{\mathcal{D}}\right)\right)}\right),
\end{equation} 
where $\mathbf{\mathcal{J}}=\mathbf{\mathcal{D}}^{\ast}\left(\Lambda_{2}-\beta_{2}\mathbf{Y}\right)+\beta_{1}\mathbf{\mathcal{X}}-\Lambda_{1}$, fftn and ifftn are 3D Fast Fourier transform and 3D inverse Fast Fourier transform, respectively. It is noticed that operator $\mathbf{\mathcal{D}}^{\ast}\mathbf{\mathcal{D}}$ can be computed outside the loop to decrease the computational cost.
\subsection{The solution of sub-problem (\ref{m})}
Considering the an-isotropic TV, the sub-problem (\ref{m}) can be the written as:
\begin{equation}\label{n}
\min_{\mathbf{Y}} \lambda \|\mathbf{Y}\|_{1} +{\frac{\beta_{2}}{2}\|\mathbf{Y}-\left(\mathbf{\mathcal{D}}\left(\mathbf{\mathcal{Z}}\right)-\frac{\Lambda_{2}}{\beta_{2}}\right)\|_{\mathrm{F}}^{2}}
\end{equation}
\begin{algorithm}[t]
	\caption{SHTRA}
	\begin{algorithmic}
		\STATE \textbf{Input:} Zero-filled observed tensor $\mathbf{\mathcal{T}}$, observed index $\mathbf{\Omega}$
		\STATE  \textbf{Initialization:} $\mathbf{\mathcal{X}}_{\Omega}=\mathbf{\mathcal{T}}_{\Omega},$ For $n = 1,\cdots, N$ random sample $\mathbf{\mathcal{G}}_{n}$, iter = 0, maxiter, residual bound $ \epsilon $, $\kappa$ = 1.01.
		\WHILE{\text{iter $<$ maxiter}}
		\STATE iter = iter + 1.
		\STATE $\mathbf{\mathcal{X}}_\text{last}=\mathbf{\mathcal{X}}$
		\STATE For $n=1, \cdots, N$, Update $\mathbf{\mathcal{G}}_{n}$, by \ref{f}
		\STATE For $n=1, \cdots, N$ Update $\mathbf{\mathcal{M}}_{n}$ by \ref{M}
		\STATE Update $\mathbf{\mathcal{Z}}$ by \ref{Z}
		\STATE Update $\mathbf{Y}$ by \ref{Y}
		\STATE Update $\mathbf{\mathcal{X}}$ by \ref{X}
		\STATE Update Dual Variables $\Lambda_{1}, \Lambda_{2}, \Lambda_{3}$ by \ref{penalty}
		\STATE $\beta=\min( \kappa \times\beta,10)$
		\STATE If $\frac{\|\mathbf{\mathcal{X}}-\mathbf{\mathcal{X}}_\text{last}\|_\text{F}}{\|\mathbf{\mathcal{X}}_\text{last}\|_\text{F}} < \epsilon,~ \textbf{Break}$
		\ENDWHILE
		\STATE \textbf{Output:} recovered tensor $\mathbf{\mathcal{X}}$
	\end{algorithmic} 
		\label{algorithm:SHTRA}
\end{algorithm}

The sub-problem (\ref{n}) can solved by soft-thresholding operator as follows:
\begin{equation}\label{Y}
\mathbf{Y}=\operatorname{sth}\left(\mathbf{\mathcal{D}}\left(\mathbf{\mathcal{Z}}\right)-\frac{\Lambda_{2}}{\beta_{2}},\frac{\lambda}{\beta_{2}}\right)
\end{equation}
where sth is the soft-thresholding operator defined as follows:
\begin{equation}
\operatorname{sth}(a,\tau)=\operatorname{sgn}\left(a\right)\max\left(|a|-\tau,0\right)
\end{equation}
where $a=\mathbf{\mathcal{D}}\left(\mathbf{\mathcal{Z}}\right)-\frac{\Lambda_{2}}{\beta_{2}}$ and $\tau=\frac{\lambda}{\beta_{2}}$.
\subsection{The solution of sub-problem (\ref{o})}
The sub-problem (\ref{o}) is a convex optimization with equality constraint. We can update $\mathbf{\mathcal{X}}$ as follows:
\begin{equation}\label{X}
X_{i_{1}, \cdots, i_{d}}=\left\{\begin{array}{ll}
{\left(\frac{\mathcal{R}_{g}+\Lambda_{1}-\beta_{1}\mathbf{\mathcal{Z}}}{1-\beta_{1}}\right)_{i_{1}, \ldots, i_{d}}} & {, i_{1}, \ldots, i_{d} \notin \mathbb{O}} \\
{\mathbf{\mathcal{T}}_{i_{1}, \ldots, i_{d}}} & {, i_{1}, \ldots, i_{d} \in \mathbb{O}}
\end{array}\right.
\end{equation}
where $\mathbf{\mathcal{R}}_g=\mathcal{F}\left(\mathbf{\mathcal{G}_1, \cdots, \mathbf{\mathcal{G}_N}}\right).$
\subsection{The solution of sub-problem (\ref{p})}
According to ADMM, the dual variables can be updated by: 
\begin{equation}\label{penalty}
\Lambda_{1}=\Lambda_{1}+\beta_{1}\left(\mathbf{\mathcal{Z}}-\mathbf{\mathcal{X}}\right)$$$$
\Lambda_{2}=\Lambda_{2}+\beta_{2}\left(\mathbf{Y}-\mathbf{\mathcal{D}}\left(\mathbf{\mathcal{Z}}\right)\right)$$$$
\Lambda_{3}^{(n)}=\Lambda_{3}^{(n)}+\beta_{3}\left(\mathbf{\mathcal{M}}_{n}-\mathbf{\mathcal{G}}_{n}\right), n=1,\cdots N.
\end{equation}
The penalty vector $ \boldsymbol{\beta}=\left[\beta_{1},\beta_{2},\beta_{3}\right]^\text{T}$ can be updated as follows \cite{cao2016total},\cite{lin2011linearized}:
$$\boldsymbol{\beta}^{(k)}=\left\{\begin{array}{ll}
{\eta_{1} \boldsymbol{\beta}^{(k-1)},} & {\text { if } \boldsymbol{\zeta}^{(k)}>\eta_{2} \boldsymbol{\zeta}^{(k-1)}} \\
{\boldsymbol{\beta}^{(k-1)},} & {\text { otherwise }}
\end{array}\right.$$
where $\boldsymbol{\zeta}^{(k)}=\left[\|\mathbf{\mathcal{Z}}-\mathbf{\mathcal{X}}\|,\|\mathbf{Y}-\mathbf{\mathcal{D}}\left(\mathbf{\mathcal{Z}}\right)\|,\sum_{n=1}^{N}\left(\|\mathbf{\mathcal{M}}_{n}-\mathbf{\mathcal{G}}_{n}\|\right)\right]^\text{T}$ in $k\text{-th}$ iteration, $\eta_{1}, \eta_{2}$ are scalar factors.

 The pseudo-codes of the SHTRA are given in Algorithm \ref{algorithm:SHTRA}.
The convergence condition is set to be $\frac{\|\mathbf{\mathcal{X}}\|_\text{F}-\|\mathbf{\mathcal{X}}_\text{last}\|_\text{F}}{\|\mathbf{\mathcal{X}}_\text{last}\|_\text{F}}< \epsilon $, where $\mathbf{\mathcal{X}}$ is the recovered tensor and $\mathcal{X}_\text{last}$ is the recovered tensor of last iteration and $ \epsilon $ is the residual bound.

\subsection{Computational Complexity}
The main computational complexity of the proposed algorithm comes from the update of $\mathcal{G}_n$. For an $N$-order tensor $\mathcal{X}^{I_1\times I_2\times \cdots \times I_N}$, the computational complexity of $\mathcal{G}_n$ is mainly from the inversion of matrix $\left(\left(\mathcal{G}_{<2>}^{(\neq n)}\right)^\text{T} (\mathcal{G}_n)_{<2>}^{(\neq n)}+ \beta_{3} \mathbf{I} \right)
$  or the matrix multiplication of $(\left(\mathcal{G}_{<2>}^{(\neq n)}\right)^\text{T} $ and $(\mathcal{G}_n)_{<2>}^{(\neq n)}$. By Assuming $I_n=I$ and $R_n=R$ for $n=1,\cdots,N$, the computational complexity of $\mathcal{G}_n$ is $O(\min(R^6, I^NR^4))$. In the proposed algorithm, $R\ll I$. Therefore, the overall complexity is $O(NTI^NR^4)$, where $T$ is the number of iterations in the proposed algorithm.

\section{Numerical Experiments}
In this section, we have conducted  several groups of experiments and compared our method with the state-of-art ones including TR-ALS \cite{wang2017efficient}, TRLRF \cite{yuan2019tensor}, STTC \cite{8421084}, TMAC-TT \cite{7859390}, SPC \cite{yokota2016smooth}, and LRTV-PDS \cite{yokota2018simultaneous}. 

To quantitatively measure the missing degree, we define the sampling ratio (SR) as follows:
$$\text{SR}=\frac{O}{\prod_{n=1}^{N}I_n}$$
 where $O$ is the number of observed entries.

The peak signal-to-noise ratio (PSNR), structural similarity index (SSIM), and relative square error (RSE) are used to measure the recovery accuracy. The PSNR is the peak signal-to-noise ratio between two images which is defined as: $10 \log_{10}\left(\frac{\text{max}^{2}}{\text{MSE}}\right)$, where $\text{max}$ is the maximum fluctuation of input image data and $\text{MSE}$ is the cumulative squared error between reconstructed and original image. $\text{MSE}$ can be defined as $\frac{\|\hat{\mathbf{\mathcal{A}}}-\mathbf{\mathcal{A}}\|_\text{F}^{2}}{N}$.The higher the PSNR, the better the quality of the recovered image. The SSIM can measure the intensity of light (i.e. luminance) and the contrast of two images. It shows how closely these features vary together between two images. The higher value for SSIM indicates better results. The RSE determines the performance by computing relative error between original tensor $\mathbf{\mathcal{A}}$ and recovered tensor $\hat{\mathbf{\mathcal{A}}}$, which can be defined as $\frac{\|\mathbf{\mathcal{A}}-\hat{\mathbf{\mathcal{A}}}\|_\text{F}}{\|\hat{\mathbf{\mathcal{A}}}\|_\text{F}}$. The lower the value of RSE, the better the performance. Besides, for the HSI, we chose MPSNR, MSSIM and   spectral angle mapper (SAM) \cite{8822607} for evaluation. Spectral angle mapper (SAM)  measures the spectral similarity by calculating the angle between the spectra $\hat{\mathbf{a}}_i$ and $\mathbf{a}_i$ over the whole spatial domain, which can be defined as 
$\text{SAM}=\frac{1}{IJ}\sum_{i=1}^{IJ}\arcos\frac{\hat{\mathbf{a}}_i^{\operatorname{T}}\mathbf{a}_i}{\|\hat{\mathbf{a}}_i\|_{\operatorname{F}}\|\mathbf{a}\|_{\operatorname{F}}}$, where $\mathcal{A}\in \mathbb{R}^{I\times J\times K}$ is the HSI and $\hat{\mathbf{a}}_i\in \mathbb{R}^{K}$ is the pixel of $\hat{\mathcal{A}}$ and $\mathbf{a}\in\mathbb{R}^{K}$ is the pixel of $\mathcal{A}$.

 All the experiments have been conducted by using MATLAB R2018b on the desktop computer with the specification Intel(R) Core(TM) i5-4590 CPU, 3301 MHz, 4 Core(s) and 8GB RAM.
Two types of datasets are used for images recovery, i.e. 1) color image dataset, 2) HSI dataset. In color image dataset, we use images with two different sizes. One group consists of standard color images with the size $256 \times 256 \times 3$ and can be seen in Fig. \ref{fig:A}., the other images are of the size $321 \times 481 \times 3$ \footnote{\url{https://www2.eecs.berkeley.edu/Research/Projects/CS/vision/bsds/}}, as can be seen in Fig. \ref{fig:B}. For HSI dataset\footnote{\url{http://www.ehu.eus/ccwintco/index.php/Hyperspectral_Remote_Sensing_Scenes}}, we use two different types. One is the Indian pines with the size $145 \times 145 \times 224$, and other is Salinas with the size $512 \times 217 \times 224$, as can be seen in Fig. \ref{2}.

\subsection{Color image with the size $256 \times 256 \times 3$ }
In this group of experiments, we apply the proposed SHTRA method along with STTC,  TRLRF, TR-ALS, SPC, LRTV-PDS and TMac-TT on four color images with the size of $256 \times 256 \times 3$. We set the trade-off parameters between smoothness and low-rank terms as 0.0005  for STTC. We set the weight of direction total variation as $\mathbf{w}=\left[4,4,0\right]^\text{T}$ for STTC in all color image experiments as provided in \cite{8421084}. We let  TR rank be 15 for TR-ALS and TRLRF in all color images based experiments. In addition, the parameter settings of SPC, TMac-TT and  LRTV-PDS follow their original papers. We randomly chose the SR from $10\%$ to $60\%$. In the proposed method, we set the trade-off parameter between total variation and low-rank term as 0.0003. We tune the penalty factors and set  $\boldsymbol{\beta} =\left[0.001; 0.001; 0.8\right]$. Besides, we also set the weight of total variation as $\mathbf{w}=\left[4; 4; 0\right]$ in all color images based experiments.  The maximum number of iterations is set to be 400 and the threshold is $ \epsilon = 5 \times 10^{-4}$. For fairness,  we also set the TR rank to be 15 for all TR-factors for the proposed method.

\subsubsection{Experimental Results  }
Fig. \ref{fig:psnrA} shows the recovery performance of 4 color images in terms of PSNR, RSE, and SSIM  with SR from 10\% to 60\%. From \ref{fig:psnrA}, it can be observed that the curves indicate the bottom-line performance of non-smooth  methods that includes TRLRF, TR-ALS, and TMAC-TT. On the other side, it can be seen that method with smoothness term such as SPC, STTC, and SHTRA show better performance compared with non-smooth methods. Therefore, it may prove that smoothness constraints can effectively enhance the recovery performance. 

In Fig. \ref{fig:psnrA}, we can also see that the proposed algorithm exhibits superior recovery performance in most cases compared to state-of-art-algorithms. Compared with other tensor-ring based methods such as TR-ALS and TRLRF,  SHTRA has achieved far better performance in every index. Fig. \ref{imageA} shows the recovered images with different sample ratios ranging from  10\% to 40\%. When the sampling ratio is very low, TRLRF has recovered the image with the features hard to identify. However, the proposed SHTRA shows better effectiveness for the reconstruction of missing entries in most cases compared with the others.

\subsection{Color image with the size $321 \times 481 \times 3$ }
In this group of experiments, we choose some color images with the size of $321 \times 481 \times 3$ for comparison of  the proposed method with TR-ALS, TRLRF, TMAC-TT, SPC, LRTV-PDS, and STTC. We set the parameters the same as those in last subsection. We also set the condition for convergence as $ \epsilon =  5 \times 10^{-4}$.

Fig. \ref{fig:psnrB} shows the recovery performance for different sample ratios from 10\% to 60\%. PSNR, RSE, and SSIM  evaluate the recovery performance. From \ref{fig:psnrB}, we can see  that as sampling ratio increases, the proposed algorithm performs better in terms of PSNR. On the other side, Fig. \ref{fig:psnrB}  indicates that the proposed algorithm has minimized the relative error between recovered image and original image as the sampling ratio increases compared with the others.

Fig. \ref{imageB} shows the recovered color images with sample ratios 30\%, 40\%, 50\%, and 60\%, respectively. From Fig. \ref{imageB}, it can be seen that as the sample ratio increases, the proposed algorithm successfully recovers the missing parts with good performance compared with the state-of-the-art ones. Besides,  the proposed SHTRA outperforms the others in terms of PSNR.


\subsection{Hyperspectral image}
In this group of experiments, we use two hyperspectral images Indianpines and Salinas with the size $145 \times 145\times 224$ and $512 \times 217\times 224$ for comparison of the proposed algorithm with STTC, SPC, LRTV-PDS, TR-ALS, and TRLRF. We tune the parameters to get the optimal performance. We set the TR rank = 10 for TR-ALS, TRLRF, and the proposed SHTRA. We set the weights for total variation as $ \mathbf{w} =[4; 4; 10]$, $ \mathbf{w} = [2; 2; 10]$ for STTC and SHTRA, respectively. We also set the trade-off parameter between low-rank term and smoothness as 0.003 and 0.0005 for STTC and SHTRA, respectively. The trade-off parameter in SHTRA is set as $\rho = 12$, and we also tune the penalty factors for the optimal performance with $\boldsymbol{\beta} = [0.001; 0.001; 0.8]$. For convergence, we the set condition $\epsilon =  10^{-4}$, and the  maximum number of iterations is set to 300.

\subsubsection{Indianpines with the size $145\times 145\times 224$}  
Indianpines with the size $145\times 145\times 224$ is chosen in this group where the spatial resolution is $145 \times 145$ and 224 refers to spectral reflectance bands.
Due to computational cost, we only chose the first 30 bands, resulting in the size  $145 \times 145\times 30$.  Table \ref{HSI} shows the quantitative evaluation of recovery with sample ratios ranging from 10\%  to 50\%. It shows that the proposed algorithm recovers the missing image Indianpines with superior performance compared to state-of-the-art ones in terms of MSSIM, MPSNR and SAM.  

Fig. \ref{HSIIndian} shows the reconstruction with the sample ratios ranging from 10\% to 40\% with different spectral band in grayscale color-map. We can see that the proposed algorithm successfully recovers the missing entries for each SR, and the performance is better than those of the others. In addition, It can be seen that when S.R=10\%, TRLRF, STTC, and LRTV-PDS fail to reconstruct the image.  However, with SR increasing, the performance of TRLRF and SPC improves along with the proposed algorithm.

\subsubsection{Salinas with the size  $512 \times 217\times 224$}
We resize Salinas by sampling its spatial resolution to $256\times 109$, and keep the first 30 bands due to computational cost. The original image is resized to $256\times 109\times 30$. 

Table \ref{HSI} shows that the proposed algorithm has achieved better performance in most cases. When SR = 20\% and SR = 40\%, TRLRF has slightly good performance compared to other algorithms. Fig. \ref{HSISaline.eps} shows the recovered  Salinas for each algorithm with different SRs. As can be seen,  LRTV-PDS and STTC have the worst performance when SR=10\%. In contrast, the recovery performance of SPC and SHTRA  is good. Besides,  SHTRA is superior to the other methods in terms of recovery resolution when  SR=50\%.

	\begin{table*}[t]
	\caption{Quantitative comparison of different algorithms for Salinas and Indianpines completion. }
	\begin{center}
			\begin{tabular}{|c|c|c|c|c|c|c|c|c|}
				\hline
				Datasets                  &     SR  & measure indexes    & STTC & SPC  & LRTV-PDS  & TR-ALS  &TRLRF & SHTRA   \\ \hline
				\multirow{20}{*}{Indianpines}  &  & M-PSNR (dB)& 39.779 &49.595& 24.528 & 47.182&42.680&\textbf{50.319}\\ \cline{3-9}
				& & M-SSIM& 0.143&0.791&0.004&0.626&0.429&\textbf{0.804}\\  \cline{3-9}
				& 10\%& SAM&0.142 &0.045&1.205&0.085 &0.124 &\textbf{0.045}\\  \cline{3-9}
				& & CPU time (sec)&52.68&58.37&1352&92.89&190.6& 444.3\\  \cline{2-9}
				&  & M-PSNR (dB)& 39.942&50.443&25.040&51.806&52.438&\textbf{ 52.558}\\ \cline{3-9}
				& & M-SSIM& 0.231&0.831&0.010&0.838&0.852&\textbf{0.880}\\  \cline{3-9}
				& 20\%& SAM& 0.141 & 0.042&1.113&0.046&0.046 &\textbf{0.038}\\  \cline{3-9}
				& & CPU time (sec)&71.09 &111.8&93.32&27.18&192.2&440.1\\  \cline{2-9}
				&  & M-PSNR (dB)& 40.332  &51.144&25.620 & 52.660 &54.046&\textbf{54.351}\\ \cline{3-9}
				& & M-SSIM&0.305&0.857&0.017&0.872&0.898&\textbf{0.916}\\  \cline{3-9}
				& 30\%& SAM & 0.137&0.04&0.995&0.041&0.039&\textbf{0.034}\\  \cline{3-9}
				& & CPU time (sec)&52.92 &97.44&60.01&20.01&177.2&310.4\\  \cline{2-9}
					&  & M-PSNR (dB)& 51.599&51.821&26.288 & 53.165 &55.132&\textbf{55.285}\\ \cline{3-9}
				& & M-SSIM&0.828&0.881&0.025&0.887&0.928&\textbf{0.929}\\  \cline{3-9}
				& 40\%& SAM & 0.053&0.038&0.888&0.038&0.031&\textbf{0.031}\\  \cline{3-9}
				& & CPU time (sec)&57.65 &83.53&36.43&20.95&35.63&122.7\\  \cline{2-9}
					&  & M-PSNR (dB)& 52.956  &52.676&27.077 & 53.361 &56.268&\textbf{56.384}\\ \cline{3-9}
				& & M-SSIM&0.867&0.902&0.034&0.893&0.946&\textbf{0.947}\\  \cline{3-9}
				& 50\%& SAM & 0.046&0.034&0.785&0.038&0.027&\textbf{0.027}\\  \cline{3-9}
				& & CPU time (sec)&53.63 &68.05&41.06&21.24&25.82&85.99\\  \cline{1-9}
				\multirow{20}{*}{Salinas}  &  & M-PSNR (dB)& 42.095 &54.399& 31.968 & 52.037&48.313&\textbf{53.615}\\ \cline{3-9}
			& & M-SSIM& 0.130&0.845&0.008&0.663&0.489&\textbf{0.79}\\  \cline{3-9}
			& 10\%& M-SAM&0.265 &\textbf{0.032}&1.209&0.101 &0.180 &0.04\\  \cline{3-9}
			& & CPU time (sec)&181.2&276.6&325.8&85.27&255.2& 561.1\\  \cline{2-9}
			&  & M-PSNR (dB)& 51.333&55.942&32.48&60.467&\textbf{61.962 }&55.790\\ \cline{3-9}
			& & M-SSIM&0.639&0.889&0.022&0.910&\textbf{0.933}&0.904\\  \cline{3-9}
			& 20\%& M-SAM&0.143& 0.029&1.113&0.025&\textbf{0.02} &0.026\\  \cline{3-9}
			& & CPU time (sec)&92.89 &85.79&53.53&15.01&100.3&281.2\\  \cline{2-9}
			&  & M-PSNR (dB)& 53.448 &56.946&33.057 & 61.523 &63.359&\textbf{63.466}\\ \cline{3-9}
			& & M-SSIM&0.710&0.910&0.039&0.935&0.952&\textbf{0.954}\\  \cline{3-9}
			& 30\%& M-SAM & 0.122&0.027&0.995&0.022&0.019&\textbf{0.018}\\  \cline{3-9}
			& & CPU time (sec)&66.95 &94.82&31.13&14.41&49.26&78.02\\  \cline{2-9}
			&  & M-PSNR (dB)& 55.268&57.911&33.727 & 61.547 &\textbf{64.686}&64.169\\ \cline{3-9}
			& & M-SSIM&0.761&0.926&0.058&0.932&\textbf{0.964}&0.958\\  \cline{3-9}
			& 40\%& M-SAM &0.104&0.026&0.887&0.021&0.017&\textbf{0.017}\\  \cline{3-9}
			& & CPU time (sec)&60.35 &87.53&32.44&214.72&35.9&64\\  \cline{2-9}
			&  & M-PSNR (dB)& 43.931  &58.874&34.52 & 61.886 &65.144&\textbf{65.631}\\ \cline{3-9}
			& & M-SSIM&0.352&0.937&0.082&0.937&0.966&\textbf{0.97}\\  \cline{3-9}
			& 50\%& M-SAM & 0.216&0.025&0.785&0.022&0.016&\textbf{0.015}\\  \cline{3-9}
			& & CPU time (sec)&34.46 &82.67&30.97&16.44&16.48&41.97\\  \cline{1-9}
			\end{tabular}
	\end{center}
	\label{HSI}
\end{table*}%

\section{Conclusion}
In this paper, we develop a low rank hierarchical tensor ring approximation for image completion. The newly proposed hierarchical tensor ring can result in more compact representation for multidimensional images in applications, in comparison with the other decompositions. By low rank approximation in both layers of the hierarchical tensor ring, automatically tuned TR ranks can alleviate overfitting when the rank is set to be large and there are a limited number of observations. To enhance the recovery performance, total variation is taken in the optimization model for exploiting local piece-wise smoothness additionally. It can be solved by ADMM.  Experimental results on color image and HSI show that the proposed algorithm outperforms the state-of-the-art ones.


%

\ifCLASSOPTIONcaptionsoff
  \newpage
\fi



%

\bibliographystyle{ieeetr}
\bibliography{bibbase}

\end{document}